\documentclass[11pt]{amsart}
\usepackage[margin=0.9in]{geometry}
\usepackage{amsmath,amssymb,amsthm,mathtools}
\usepackage{stmaryrd}
\usepackage{booktabs,array}
\usepackage{graphicx}
\usepackage{xcolor}
\usepackage[numbers]{natbib}
\usepackage[colorlinks=true,linkcolor=blue!50!black,citecolor=blue!50!black,urlcolor=blue!50!black]{hyperref}

\newtheorem{theorem}{Theorem}[section]
\newtheorem{proposition}[theorem]{Proposition}
\newtheorem{lemma}[theorem]{Lemma}
\newtheorem{corollary}[theorem]{Corollary}

\theoremstyle{definition}

\theoremstyle{remark}
\newtheorem{remark}[theorem]{Remark}

\newcommand{\R}{\mathbb{R}}
\newcommand{\Prob}{\mathbf{P}}
\newcommand{\E}{\mathbf{E}}
\newcommand{\obj}{\mathrm{obj}}
\newcommand{\val}{\mathrm{val}}
\newcommand{\cP}{\mathcal{P}}
\newcommand{\cR}{\mathcal{R}}

\title{Flag algebras and the exact largest probability of a linear inequality between iid random variables}

\author{Anton Zaytsev}
\email{anton.a.zaytsev@gmail.com}
\author{Pierre C. Bellec} 

\begin{document}

\begin{abstract}
How large can $\Prob[X_1+X_2+X_3<2X_4]$ be when $X_1,\dots,X_4$ are independent with a common law $\mu$ on $[0,\infty)$?
Previous work gave lower bounds through explicit constructions and upper bounds through exchangeability.
The best lower bound, $0.4006953723\ldots$, was conjectured
by \citet{BF} to be the exact value of the supremum.
We prove this conjecture and formalize the proof in Lean~4.
The improved upper bound uses the fact that
independence gives more than exchangeability: two subsamples are conditionally independent given their shared points. As in Razborov's flag algebras, this yields a semidefinite program on the sign patterns of a few sample points. The new ingredient is to record also whether one sample point exceeds a fixed multiple of another, at ratios read off the extremal construction: with three such cutoffs, four sample points give an exact certificate, while additive comparisons alone leave a gap at seven points. We extend this study to the functions
\[
S_2(c) = \sup_\mu \Prob_\mu[X_1+X_2<cX_3],\qquad
S_3(c) = \sup_\mu \Prob_\mu[X_1+X_2+X_3<cX_4]
\]
and determine their exact values on several intervals of $c$.
On other intervals, we provide lower and upper bounds that
certify eighteen discontinuity points of $S_2$ and five of $S_3$,
consistent with a nested staircase structure.
The certified jumps are explained by new winning patterns becoming feasible
in the lower-bound constructions.
We answer Problem~5.3 of \cite{BF} by showing that no coefficient $c$
satisfies $S_3(c)=1/2$, and prove the analogous statement for $S_2$.
The thresholds $\sup\{c:S_2(c)\le1/2\}$ and
$\sup\{c:S_3(c)\le1/2\}$ are respectively the golden ratio
$(1+\sqrt5)/2$ and the silver ratio $1+\sqrt2$.
\end{abstract}

\maketitle

\setcounter{tocdepth}{2}
%\vspace{-0.5in}
{\small
\tableofcontents
}
\clearpage

%======================================================================
\section{Introduction and main results}
\addtocontents{toc}{\protect\setcounter{tocdepth}{1}}
%======================================================================

\subsection*{Three-against-one (Figure~\ref{fig:three-against-one})} In July 2024, \citet{MO474916} asked on MathOverflow
for the value of 
\begin{equation}
    \label{C}
  C \;=\; \sup_\mu \Prob_{\mu^{\otimes 4}}\bigl[X_1+X_2+X_3<2X_4\bigr],
\end{equation}
with the supremum over all laws $\mu$ on $\R_+=[0,\infty)$, with $X_1,\dots,X_4$ independent with law $\mu$.
This is the best possible chance that one roll of the common die is more than one and a half times the average of three other rolls. We call this the four-roll game. 
Prior work established lower and upper bounds:
\cite{BF} proved
\begin{equation}\label{eq:BF}
  C_{\mathrm{BF}} \;\le\; C \;\le\; \tfrac{673}{1615} < 0.417,
\end{equation}
where
\begin{equation}\label{eq:lower}
  C_{\mathrm{BF}} \;=\; \sup_{p\in(0,1)} \frac{p(2-p)}{1+p+p^2+p^3} \;=\; 0.4006953723\ldots
\end{equation}
Conjecture~5.1 in \cite{BF} asserts that $C=C_{\mathrm{BF}}$.
This paper proves that conjecture and formalises it in Lean~4.

Independently, Tomar \cite{Tomar} proved $C\le0.40069634$ with a nine-point certificate of the kind described in Section~\ref{sec:method-flag-standard}, with a Lean~4 proof.
Section~\ref{sec:four-exact} proves the exact value $C=C_{\mathrm{BF}}$
with a four-point flag certificate that also records three pairwise ratio cutoffs,
a technique we will refer to as ``colored certificates''.

This lower bound $C_{\mathrm{BF}}$ comes from a self-similar multiscale construction \cite{BF,Yu}. It starts with the mixture $p\delta_0+(1-p)U$, where $p\in(0,1)$, $\delta_0$ is the point mass at $0$, and $U$ is the uniform law on $\{1-2^{-k}:k=1,\ldots,N\}$. This pattern is repeated at infinitely many smaller scales by taking
$X=\sum_{k=0}^{\infty}\eta^kY_k$
where the $Y_k$ are iid with law $p\delta_0+(1-p)U$ and $\eta\in(0,1)$. Taking $\eta\to0$ and $N\to\infty$ gives the lower bound. The first version of \cite{BF} considered only $p=1/2$, giving $0.4$; \cite{Yu} observed that optimizing over $p$ improves this by $0.0006953723\ldots$, yielding $C_{\mathrm{BF}}$.

The upper bound in \cite{BF} is the worst case over a single sample: for $m$ points, the largest fraction of the instances of the inequality that can hold simultaneously. Probabilistically, this is the best bound that uses only the exchangeability of $m$ draws \cite{Lacker}; \citet{BF} computed it for $m=20$ by mixed integer programming. The problem is also Problem~6.39 of the AlphaEvolve collection \cite{GGTW}, where automated search produced lower bounds only. In a recent talk, \citet{Fritz26} conjectured that ``an exact computation method based on the fractal nature of near-optimal distributions exists''.

We also vary the right-hand side coefficient and call the function
$$
  S_3(c) \;=\; \sup_\mu \Prob_{\mu^{\otimes 4}}\bigl[X_1+X_2+X_3<cX_4\bigr].
$$
the ``three against one'' game. See Figure~\ref{fig:three-against-one}
for a summary of the bounds we obtain on this function.
Section~\ref{toplevel_sec:thmD} develops constraints from near-optimal
laws and uses them to prove five additional exact $S_2$ plateaus.
Four-point colored certificates also prove the first $S_3$
plateau above $2$ exactly, through $(1+\sqrt{13})/2$, and an exact interval
immediately after that jump, through $2.35$
(Section~\ref{sec:three}). Below $2$, a smaller certificate proves
an exact plateau on $(1,\sqrt2]$.

\begin{figure}[b]
\centering
\includegraphics[width=\textwidth]{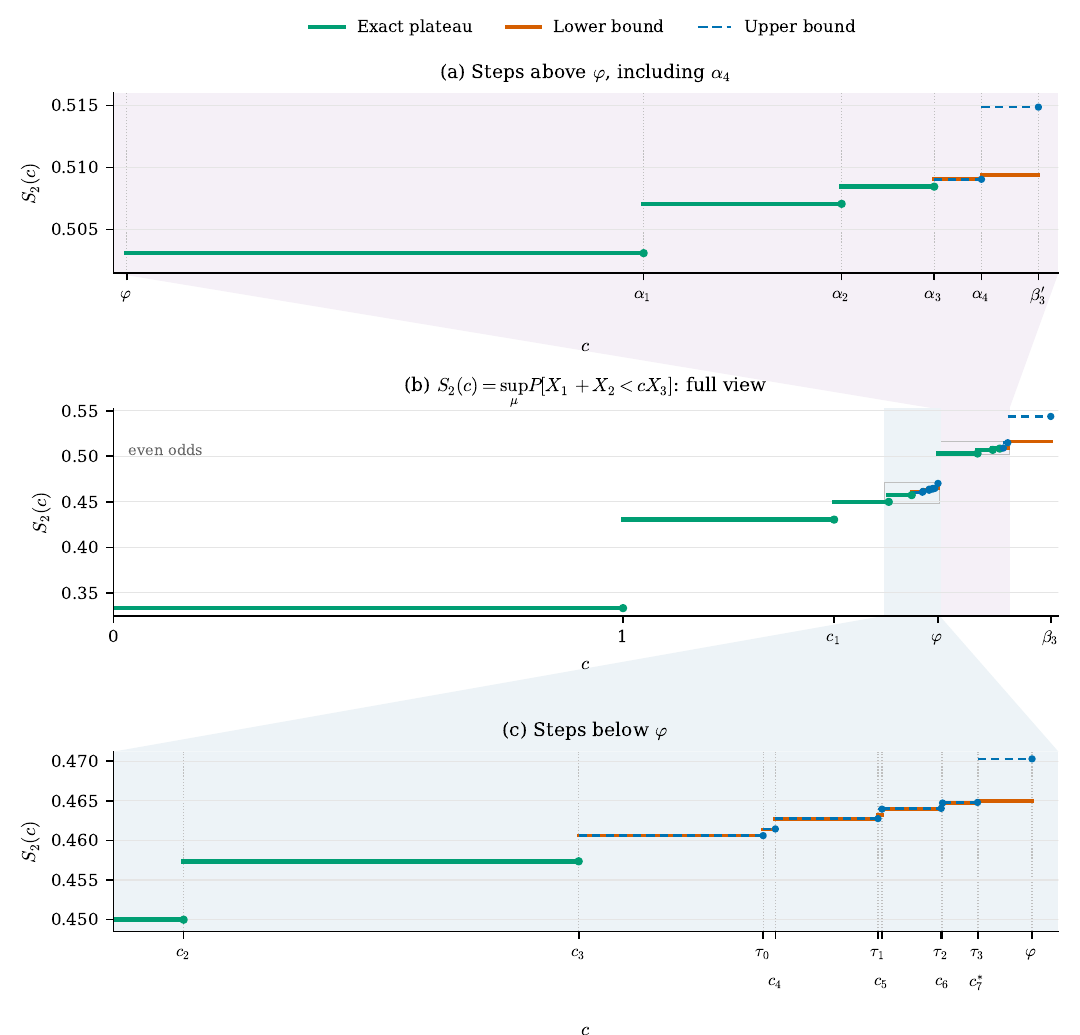}
\caption{Bounds for $S_2$ from Table~\ref{tab:steps-bounds}: exact plateaus (green), lower bounds (orange), and upper bounds (blue, dashed). Each exact plateau is drawn once. Shaded connectors link the middle full view to the top zoom above $\varphi$ through $\alpha_4$ and the bottom zoom below $\varphi$. The guide $c_7^*$ marks the next ladder threshold within $(\tau_3,\varphi]$, using the bounds for that whole interval.}
\label{fig:staircase}
\end{figure}

\subsection*{Two-against-one (Figure~\ref{fig:staircase})} 
Similarly to $S_3(c)$, we study the ``two-against-one'' function
\begin{equation}
    \label{S2}
  S_2(c) \;=\; \sup_\mu \Prob_{\mu^{\otimes 3}}\bigl[X_1+X_2<cX_3\bigr].
\end{equation}
Figure~\ref{fig:staircase} summarizes the exact plateaus and lower/upper bounds
we obtain on this function.
The upper bound on the first plateau, for $c\in(0,1]$, follows
from exchangeability:
$\mathbf 1\{X_1+X_2<c X_3\}+\mathbf 1\{X_1+X_3<c X_2\}+\mathbf 1\{X_2+X_3<cX_1\}
\le 1$ for $c\le1$, so the average of the three terms is at most $1/3$
and $S_2(c)\le1/3$. The idea of flag algebras is to capture more
subtle structure of an iid sample by conditioning on a subset.
%A simpler problem is to establish
%$\inf_\mu \Prob_{\mu^{\otimes 3}}\bigl[2X_3\ge X_1+X_2\bigr] = 1/3$,
%which has the following casino game interpretation:
%Three cups each hold a die; the dice are identical but possibly nonstandard. You pick a cup and win if your die shows at least the average of the other two. Is the game fair? The casino can make it unfavourable: with suitable dice your chance of winning comes arbitrarily close to $1/3$, and no closer \cite{BF}. 

\subsection*{The idea}
Independent draws have one more property than exchangeable ones. Fix a set $R$ of points of an iid sample, the roots, and two further sets $A$ and $B$ of the same size, disjoint from each other and from $R$. What we look at on a set of points is its sign pattern: which instances of the linear inequalities of the problem hold among these points. Given $R$, the pattern on $R\cup A$ is independent of the pattern on $R\cup B$, and the two have the same conditional law. So for any score $v$ of a pattern,
\[
  \E\bigl[v(R\cup A)\,v(R\cup B)\bigr]\;=\;\E\Bigl[\bigl(\E[v(R\cup A)\mid R]\bigr)^2\Bigr]\;\ge\;0 .
\]
These inequalities hold for iid samples and, by de Finetti's theorem, for mixtures of them, but not for all exchangeable samples. They are the positive semidefinite constraints of Razborov's flag algebras \cite{Razborov07,Razborov13}. Section~\ref{sec:example} carries out the whole argument on four points, by hand.

We are free to choose which relations a pattern records, and the choice turns out to matter. Section~\ref{sec:pair-colors} also records pairwise ratio colors, that is, comparisons $x_j<qx_i$ for a few fixed $q$. This gives an exact three-point certificate for the plateau of $S_2$ on $(1,\sqrt2\,]$ and, with three cutoffs on four points, the exact value of $C$ (Section~\ref{sec:four-exact}).

\subsection{Results}

\begin{theorem}[The four-roll game]\label{thm:main}
$C=C_{\mathrm{BF}}$.
\end{theorem}

Recording only additive comparisons, the flag algebra SDP gives $C<0.4025452$ with six points and $C<0.4007545$ with seven, leaving a gap of less than $5.92\cdot10^{-5}$ above the construction of \citet{BF} (Section~\ref{sec:bta}). Recording also the pairwise ratio cutoffs $7/6$, $3/2$, and $2$ closes the gap with just four points. The certificate consists of four positive semidefinite matrices over the algebraic field defining $C_{\mathrm{BF}}$, and bounds the success probability for every law. Section~\ref{sec:four-exact} gives the short proof, formalized in Lean~4. Section~\ref{sec:thmD} develops additional constraints from near-optimal laws and a fixed-point inequality $S\le\psi(S)$ for the unknown optimum $S$. These tools prove five further exact two-against-one plateaus in Table~\ref{tab:steps-bounds}.

The method applies to any event cut out by strict linear inequalities, and we use it on two one-parameter families from \cite[(3.4)]{BF}. For $1<c<2$ consider the two-against-one game $S_2(c)$ in
\eqref{S2}.
In casino terms: the house rolls once, you roll twice, and the house wins if $c$ times its roll exceeds your total; the house chooses the (common) dice, and $S_2(c)$ is its best chance. By analogy with \cite[Problem~5.3]{BF}, call $c$ \emph{fair} if this best chance $S_2(c)$ is exactly even, $1/2$. Fairness concerns the best law: particular dice can always give exactly $1/2$ once some law gives more, since mixing that law with a point mass at $0$ moves the probability continuously down to $0$. V.~Yu observed that the analogous function for three against one is discontinuous, so that $1/2$ need not be attained, and suggested studying the threshold $\sup\{c:S(c)\le1/2\}$ instead \cite{Yu}. For two against one this is exactly what happens: up to the golden ratio even the house's best dice leave it below even, and just beyond, it can do better than even. Let $\varphi=(1+\sqrt5)/2$.

\begin{theorem}[No fair coefficient for two against one]\label{cor:fair2}
Put
$L_\varphi=\frac29\big(\sqrt{75-30\sqrt3}+2\sqrt3-6\big)
           =0.5030984086353\ldots.$
Then
\[
 S_2(c)\le0.4702903<\tfrac12\quad(0<c\le\varphi),\qquad
 S_2(c)\ge L_\varphi>\tfrac12\quad(c>\varphi).
\]
Let $\alpha_1=1.6956208\ldots$ be the unique positive root of
$c^3=c^2+2$. We have $S_2(c)=L_\varphi$ for
$\varphi<c\le\alpha_1$, so
$S_2(\varphi+)=L_\varphi$ and
$S_2(\varphi+)-S_2(\varphi)>0.0328081$.
Consequently
$\sup\{c>0:S_2(c)\le\tfrac12\}=\varphi$,
and no positive coefficient satisfies $S_2(c)=1/2$.
\end{theorem}

The upper bound in Table~\ref{tab:steps-bounds} at $c=\varphi$ and the uniform right-side
construction are certified at the exact golden-ratio endpoint. The golden
ratio appears because it is where the stronger comparison $r+1<c$
first becomes compatible with $r>1/c$ in the two-atom construction (Section~\ref{sec:windows} and
Figure~\ref{fig:windows}). Below $\varphi$, $S_2$ appears to climb a
staircase (Figure~\ref{fig:staircase}).

\begin{theorem}[Two-against-one]\label{thm:steps}
The function $S_2(c)$ has right-hand discontinuities at
$$
\beta_1=1,\quad c_1,\ldots,c_6,\quad \tau_0,\ldots,\tau_3,\quad
\beta_2=\varphi,\quad \alpha_1,\alpha_2,\alpha_3,\alpha_4,\quad \beta_3',\beta_3,
$$
defined in Table~\ref{tab:steps-thresholds}, with bounds on $S_2(c)$ and
on the discontinuity jumps given in Table~\ref{tab:steps-bounds}.
\end{theorem}

\begin{table}[t]
\centering
\small
\setlength{\tabcolsep}{2.5pt}
\begin{tabular}{cc|l|c}
\toprule
threshold & sub-threshold & definition & approximate value\\
\midrule
$\beta_1=1$ & & positive root of $c^2=1$ & $1$\\
           &$c_1=\sqrt2$ & positive root of $c^2=2$ & $1.4142136$\\
           &$c_2$ & positive root of $c^3=c+2$ & $1.5213797$\\
           &$c_3$ & positive root of $c^4=c^2+c+2$ & $1.5663833$\\
           &$\tau_0=\sqrt[3]4$ & positive root of $c^3=4$ & $1.5874011$\\
           &$c_4$ & positive root of $c^5=c^3+c^2+c+2$ & $1.5887905$\\
           &$\tau_1$ & positive root of $c^4=c^2+4$ & $1.6004852$\\
           &$c_5$ & positive root of $c^6=c^4+c^3+c^2+c+2$ & $1.6009326$\\
           &$\tau_2$ & positive root of $c^5=c^3+c^2+4$ & $1.6076825$\\
           &$c_6$ & positive root of $c^7=c^5+c^4+c^3+c^2+c+2$ & $1.6078359$\\
           &$\tau_3$ & positive root of $c^6=c^4+c^3+c^2+4$ & $1.6118225$\\
$\beta_2=\varphi=(1+\sqrt5)/2$ & & positive root of $c^2=c+1$ & $1.6180340$\\
                       &$\alpha_1$ & positive root of $c^3=c^2+2$ & $1.6956208$\\
                       &$\alpha_2$ & positive root of $c^4=c^3+c+2$ & $1.7253402$\\
                       &$\alpha_3$ & positive root of $c^5=c^4+c^2+c+2$ & $1.7392343$\\
                       &$\alpha_4$ & positive root of $c^6=c^5+c^3+c^2+c+2$ & $1.7463315$\\
$\beta_3'$ & & positive root of $c(c-1)^2=1$ & $1.7548777$\\
$\beta_3$ & & positive root of $c^3=c^2+c+1$ & $1.8392868$\\
\bottomrule
\end{tabular}
\caption{The eighteen discontinuity points in Theorem~\ref{thm:steps}, listed in increasing order.}
\label{tab:steps-thresholds}
%\end{table}

%\begin{table}[b]
\centering
\small
\setlength{\tabcolsep}{3pt}
\begin{tabular}{c|c|c|c}
\toprule
interval for $c$ & certified bounds for $S_2(c)$ & gap & discontinuity\\
\midrule
$(0,1\,]$
  & $S_2(c)=\tfrac13$ & $0$
  & $S_2(1+)-S_2(1)>9.7\mathrm{e}{-2}$\\
$(1,c_1\,]$
  & $S_2(c)=(2\sqrt7-4)/3$ & $0$
  & $S_2(c_1+)-S_2(c_1)>1.9\mathrm{e}{-2}$\\
$(c_1,c_2\,]$
  & $S_2(c)=0.44998962221368\ldots$ & $0$
  & $S_2(c_2+)-S_2(c_2)>7.3\mathrm{e}{-3}$\\
$(c_2,c_3\,]$
  & $S_2(c)=0.45736151856620\ldots$ & $0$
  & $S_2(c_3+)-S_2(c_3)>3.2\mathrm{e}{-3}$\\
$(c_3,\tau_0\,]$
  & $0.4606036\le S_2(c)\le0.4606037$ & $4.6\mathrm{e}{-9}$
  & $S_2(\tau_0+)-S_2(\tau_0)>7.8\mathrm{e}{-4}$\\
$(\tau_0,c_4\,]$
  & $0.4613898\le S_2(c)\le0.4614455$ & $5.6\mathrm{e}{-5}$
  & $S_2(c_4+)-S_2(c_4)>1.3\mathrm{e}{-3}$\\
$(c_4,\tau_1\,]$
  & $0.4627479\le S_2(c)\le0.4627493$ & $1.2\mathrm{e}{-6}$
  & $S_2(\tau_1+)-S_2(\tau_1)>4.7\mathrm{e}{-4}$\\
$(\tau_1,c_5\,]$
  & $0.4632210\le S_2(c)\le0.4639491$ & $7.3\mathrm{e}{-4}$
  & $S_2(c_5+)-S_2(c_5)>1.3\mathrm{e}{-6}$\\
$(c_5,\tau_2\,]$
  & $0.4639503\le S_2(c)\le0.4640382$ & $8.8\mathrm{e}{-5}$
  & $S_2(\tau_2+)-S_2(\tau_2)>2.2\mathrm{e}{-4}$\\
$(\tau_2,c_6\,]$
  & $0.4642668\le S_2(c)\le0.4647187$ & $4.5\mathrm{e}{-4}$
  & $S_2(c_6+)-S_2(c_6)>8.1\mathrm{e}{-7}$\\
$(c_6,\tau_3\,]$
  & $0.4647195\le S_2(c)\le0.4647946$ & $7.5\mathrm{e}{-5}$
  & $S_2(\tau_3+)-S_2(\tau_3)>1.5\mathrm{e}{-4}$\\
$(\tau_3,\varphi\,]$
  & $0.4649467\le S_2(c)\le0.4702903$ & $5.3\mathrm{e}{-3}$
  & $S_2(\varphi+)-S_2(\varphi)>3.2\mathrm{e}{-2}$\\
$(\varphi,\alpha_1\,]$
  & $S_2(c)=0.50309840863532\ldots$ & $0$
  & $S_2(\alpha_1+)-S_2(\alpha_1)>3.9\mathrm{e}{-3}$\\
$(\alpha_1,\alpha_2\,]$
  & $S_2(c)=0.50707388333268\ldots$ & $0$
  & $S_2(\alpha_2+)-S_2(\alpha_2)>1.3\mathrm{e}{-3}$\\
$(\alpha_2,\alpha_3\,]$
  & $S_2(c)=0.50846036134736\ldots$ & $0$
  & $S_2(\alpha_3+)-S_2(\alpha_3)>5.9\mathrm{e}{-4}$\\
$(\alpha_3,\alpha_4\,]$
  & $0.5090541\le S_2(c)\le0.5090544$ & $2.3\mathrm{e}{-7}$
  & $S_2(\alpha_4+)-S_2(\alpha_4)>3.2\mathrm{e}{-4}$\\
$(\alpha_4,\beta_3'\,]$
  & $0.5093764\le S_2(c)\le0.5148734$ & $5.5\mathrm{e}{-3}$
  & $S_2(\beta_3'+)-S_2(\beta_3')>1.3\mathrm{e}{-3}$\\
$(\beta_3',\beta_3\,]$
  & $0.5162682\le S_2(c)\le0.5437506$ & $2.7\mathrm{e}{-2}$
  & $S_2(\beta_3+)-S_2(\beta_3)>4.0\mathrm{e}{-4}$\\
\bottomrule
\end{tabular}
\caption{Certified bounds and discontinuity jumps for Theorem~\ref{thm:steps}. A zero gap denotes an exact plateau. Positive gaps are computed before display rounding and rounded to two significant digits; the unrounded gap on $(c_3,\tau_0]$ is less than $4.63\cdot10^{-9}$.}
\label{tab:steps-bounds}
\end{table}

The golden ratio $\beta_2=\varphi$ is also $\beta_2=\lim_{n\to\infty}c_n$, where $c_n$ is the positive root of $c^{n+1}=c^{n-1}+c^{n-2}+\cdots+c+2$. In the inequality rows of Table~\ref{tab:steps-bounds}, lower bounds are rounded down and upper bounds up; the equality rows display exact formulas or initial decimal digits of exact algebraic values. Each interval $(a,b]$ uses the uniform lower construction valid for every $c>a$ and the upper certificate at $b$, extended to the interval by monotonicity. The jump bounds are computed from the unrounded certificates. The jump at $\beta_3$ uses the construction giving $S_2(c)\ge0.5441582$ for every $c>\beta_3$.

The five additional zero-gap rows have exact values
\[
 V(\beta)=\frac23\left(\sqrt{7-15\beta+9\beta^2}-2+3\beta\right),
\]
where, in order of appearance, the within-scale success probabilities are
\[
 \frac4{27},\qquad \frac{\sqrt3}{9},\qquad
 \frac{2\sqrt3}{9},\qquad \frac{216+38\sqrt{19}}{961},\qquad
 2s-3s^2+(1+\gamma)s^3,
\]
with $\gamma=\sqrt3/9$ and
$s=(3-\sqrt{3-6\gamma})/(3(1+\gamma))$ in the last case.
These are the recursive lower values of Proposition~\ref{prop:onescale}.
Matching four-point pair-color certificates use the near-optimality
constraints and fixed-point argument developed in Section~\ref{sec:thmD}.
Section~\ref{sec:near-plateaus} gives their application to these five intervals,
with exact equality at $V(\beta)$ and a strictly positive margin at
$2/3$. Their covers, algebraic identities, scalar multipliers, and
positive semidefinite matrices are checked by the standalone scripts
\path{verify_second_plateau.py}, \path{verify_third_plateau.py},
\path{verify_alpha1_plateau.py}, \path{verify_alpha2_plateau.py}, and
\path{verify_alpha3_plateau.py} in \path{anc/proofs/s2/}.
The narrow bands at $\tau_0$ and $\tau_1$ use the full-strength
substitution and degree constraints; \path{verify_near_small.py}
checks their rational certificates exactly. A nonzero gap in these
rows does not establish plateau equality.
Likewise, \path{verify_alpha4_band.py} checks the band on
$(\alpha_3,\alpha_4]$, whose unrounded width is less than
$2.33\cdot10^{-7}$. The uniform right-hand construction at $\alpha_4$
gives $S_2(c)>0.5093764319755$ for every $c>\alpha_4$;
\path{verify_alpha4_discontinuity.py} checks its separation from the
endpoint upper bound, proving a jump greater than $0.000322066961$.

The exact value in the second row of Table~\ref{tab:steps-bounds} is proved in Section~\ref{sec:plateau}: recording one pairwise ratio cutoff gives an explicit optimal rank-one flag certificate on three points. The additional endpoint bounds use four-point pair-color certificates, as described in Section~\ref{sec:colored-jumps}. The thresholds $c_1$, $c_2$ and $\varphi$ are the roots of $c^2=2$, $c^3=c+2$ and $c^2=c+1$, the points where new windows for the within-scale ratios open.

For the three-against-one game $S_3(c)=\sup_\mu\Prob[X_1+X_2+X_3<cX_4]$, $2<c<3$, \citet[Problem~5.3]{BF} ask whether some $c^*$ makes it fair. Yu observed that $S_3(c)\ge L^*=0.4827817\ldots$ for all $c>2$ and asked whether this is the true value of $\lim_{c\to2^+}S_3(c)$ \cite{Yu}.

\begin{theorem}[Three against one]\label{thm:fair3}
$S_3(c)=1/4$ for $0<c\le1$.
For $1<c\le\sqrt2$, its value is $L_0=0.3096847954174\ldots$,
where $L_0$ is defined in \eqref{eq:s3-below-two}.
There is a jump greater than $0.011685187$ immediately after $\sqrt2$.
$S_3(c)=L^*=0.4827817466286\ldots$ for every
$c\in(2,(1+\sqrt{13})/2]$. In particular,
$\lim_{c\to2^+}S_3(c)=L^*$, and there is a jump strictly to the right
of $(1+\sqrt{13})/2$. We also have
$S_3(c)=L_2=0.4911560818137\ldots$ for
$(1+\sqrt{13})/2<c\le47/20$, where $L_2$ is the two-positive-atom
value in \eqref{eq:s3-second-value}. At the silver ratio,
\[
  0.4952063\;<\;S_3(1+\sqrt2)\;<\;0.4952221\;<\;\tfrac12,
\]
while $S_3(c)>0.503419372$ for every $c>1+\sqrt2$. Hence no $c$ makes $S_3(c)$ exactly $1/2$.
\end{theorem}

\begin{figure}[t]
\centering
\includegraphics[width=\textwidth]{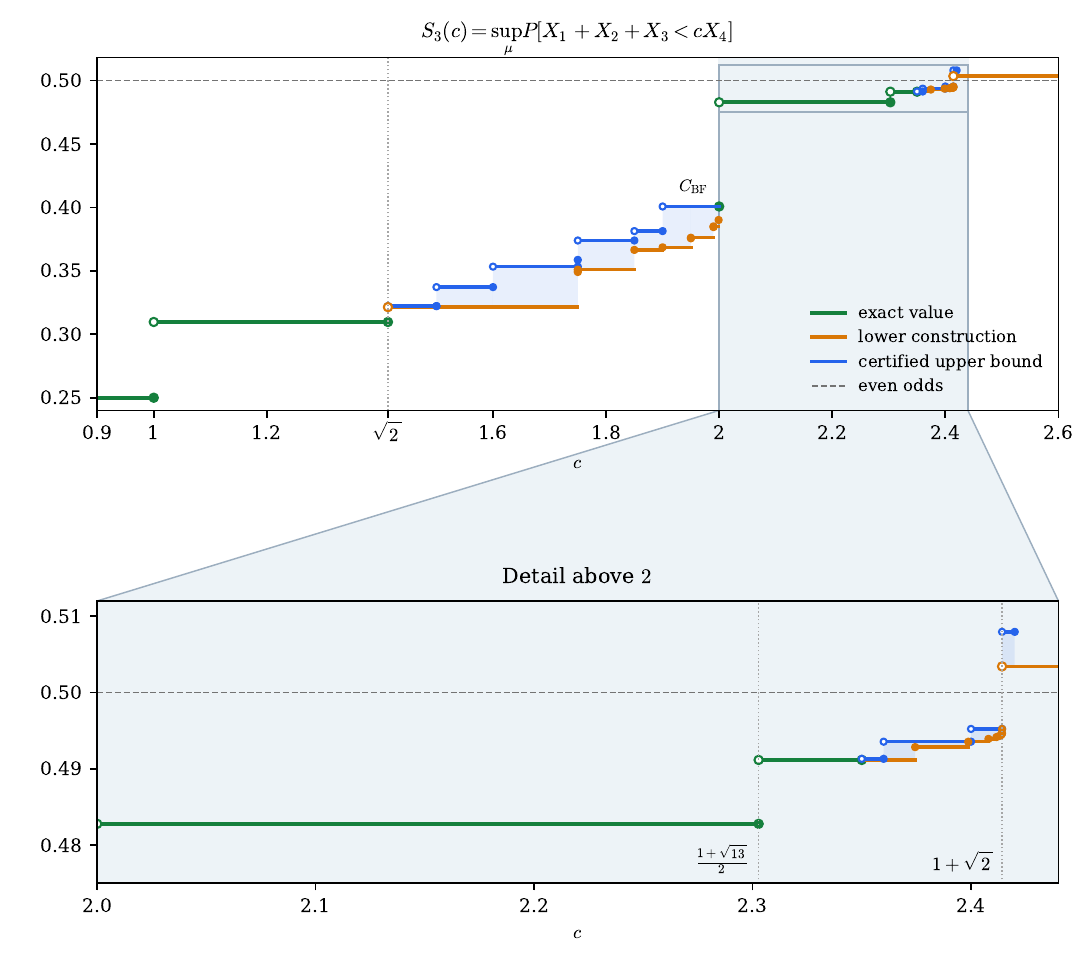}
\caption{Exact values (green), certified upper bounds (blue), and lower
constructions (orange) for $S_3$, with the remaining uncertainty shaded.
The full view above starts at $c=0.9$; the lower panel magnifies
$2\le c\le2.44$. Matching shaded backgrounds and diagonal connectors
identify the magnified region.
The exact plateaus are $(0,1]$, $(1,\sqrt2]$,
$(2,(1+\sqrt{13})/2]$, and $((1+\sqrt{13})/2,2.35]$.
The green point at $2$ is $C_{\mathrm{BF}}$; a plateau immediately to its
left remains unresolved. The plateau through $2.35$ may extend farther.
Other steps are bounds, including the ladder lower constructions below $2$.
The upper bound at $1+\sqrt2$ is below $1/2$, while the uniform lower
bound immediately to its right is above $1/2$, excluding a fair
coefficient. Open circles mark excluded endpoints.}
\label{fig:three-against-one}
\end{figure}

This answers Problem~5.3 of \cite{BF} and Yu's question about the right
limit at $2$. All three nontrivial exact intervals use ordinary four-point
pair-color certificates; the maximal right endpoint of the plateau
at $L_2$ remains unknown. The threshold $1+\sqrt2$ plays the role of $\varphi$: both are where the same kind of window opens (Figure~\ref{fig:windows}). The upper bound at $1+\sqrt2$ is a seven-point certificate whose completeness certificates are checked in $\mathbb Q(\sqrt2)$, and the lower bound is a ladder of atoms below $\sqrt2$ repeated on separated scales (Section~\ref{sec:three}); the two differ by $1.6\cdot10^{-5}$.

The certificates also give quantitative restrictions on the sample patterns of near-optimal laws: a law approaching a sharp bound must concentrate on patterns with small certificate slack (Remark~\ref{sec:structure}). The separated-scale constructions of \cite[Section~5]{BF} can guide the choice of relations and matrix kernels, as suggested by \citet{Fritz26}.

\subsection*{Related work}
\emph{De Finetti and sums of squares.} The level of our hierarchy without semidefinite constraints is exactly the bound of \citet{BF}. It goes back to Lacker's sampling-without-replacement argument \cite{Lacker} and converges by \cite[Theorem~2.4]{BF}. Relaxing a product measure by an exchangeable one and then applying sums of squares to the resulting finite problem is the approach of \cite{LevinChandrasekaran} for polynomial objectives. There it also reproduces flag-algebra bounds for graph densities. A recent paper uses a positivity condition without a shared root set \cite{TranSoh}. Related ideas drive the inflation technique of causal inference \cite{WSF19}.

\emph{Flag algebras of iid points.} Flag algebras have been applied before to structures induced by iid points, such as planar order types \cite{Goaoc} and permutations. The packing density of a permutation pattern is the supremum, over laws on the square, of the probability that iid points form that pattern \cite{PS10}, and flag algebras bound it \cite{SS18}. Flag algebras of permutations have also settled other extremal questions \cite{BHLPUV15}. In these settings the events are pure order relations, and for permutations every pattern is realizable.

\emph{What is new here.} Two features set these problems apart from earlier flag-algebra applications. First, the linear forms $x_i+x_j-cx_l$ carry additive structure, so which sign patterns can occur is itself a question of linear feasibility. We settle it with explicit Gordan certificates, which makes the list of patterns part of the checkable proof. Second, near-optimal laws live on many scales. Our semidefinite constraints are the rooted ones of flag algebras, and the rounding to exact certificates is standard \cite{PeyrlParrilo,FRV13}.

Tomar's nine-point certificate \cite{Tomar}, mentioned in the introduction, uses a closely related relaxation. Tomar asks whether a certificate of this kind reaches the exact value at some finite sample size \cite[Section~7]{Tomar}. Section~\ref{sec:four-exact} shows that four points suffice once pairwise ratio cutoffs are also recorded, using flag positivity only; for additive comparisons alone the question remains open. The near-optimality constraints of Section~\ref{sec:thmD} prove further exact two-against-one plateaus. Copositivity was used for iid inequalities of the ``123 theorem'' type in \cite{KMM14}.

%======================================================================
\addtocontents{toc}{\protect\setcounter{tocdepth}{2}}
\section{AI usage and verification: Lower bounds, upper bounds, and Lean}
\label{sec:ai-usage}

Under human guidance and supervision, frontier AI models have been used extensively to derive the results of this paper, including theorem statements, proofs, and numerical code for finding useful SDP certificates. These tools allowed us to make progress on these probability problems without prior expertise in flag algebras.
Credit is due to A. Razborov and the other mathematicians who developed and popularized flag algebras, and to the researchers responsible for these models.

For the problems considered here, sharper lower and upper bounds on
suprema over probability measures such as \eqref{C} admit explicit witnesses
that can be checked independently of the search that produced them.
\begin{itemize}
    \item Lower bounds require explicit constructions of probability distributions, typically discrete, for which the probabilities of the events of interest can be computed exactly. For constructions involving infinitely many scales, this also requires justifying the summation or recursion used to evaluate the probability.
    \item Upper bounds are obtained through semidefinite programming relaxations. Verifying a certificate requires checking that the pattern list covers every realizable case, recomputing the flag coefficients, and checking positive semidefiniteness and the certificate inequalities in exact arithmetic. These finite checks rely on the mathematical argument relating the relaxation to the original probability problem. For the exact value of \eqref{C}, Section~\ref{sec:four-exact} uses the ordinary flag-positivity argument. The additional exact $S_2$ plateaus also use the near-optimality constraints and fixed-point argument of Section~\ref{sec:thmD}.
\end{itemize}

We provide explicit lower-bound constructions and upper-bound certificates to allow independent verification. Our confidence in the AI-assisted results rests on these checks and the accompanying proofs, rather than on the reliability of the models that proposed them.

For the problem of determining the constant \eqref{C}, originally posed by \citet{MO474916}, we have also formalized the proof of Theorem~\ref{thm:main} in Lean~4. This formalization covers the probability-theoretic argument as well as the finite certificate checks. The direct four-point proof of Section~\ref{sec:four-exact} uses kernel-only checks. The Lean formalization concerns Theorem~\ref{thm:main}; the other bounds in the paper are verified by the mathematical proofs and certificate checks described above.

%======================================================================
\section{The method of flag algebras}
\label{sec:method-flag-standard}

The method turns the statement ``the sample is iid'' into finitely many inequalities that a computer, or for small cases a person, can check. We first carry it out on an example small enough to do by hand (Section~\ref{sec:example}), then state it in general (Section~\ref{sec:general}) and compare it with the bound from exchangeability (Section~\ref{sec:relation}). The last three subsections apply it to the four-roll game with six and seven points.

\subsection{A small example}\label{sec:example}
%======================================================================

We prove $S_2(6/5)\le 54/121$ with four sample points. The bound of \citet{BF} with four points gives only $1/2$. Throughout this subsection $c=6/5$. The argument has four steps: we see which point of a triple can win, write down one square whose mean is nonnegative because the sample is iid, list the fourteen patterns that four points can form, and add the square to the objective pattern by pattern, so that the worst pattern bounds the probability.

\subsubsection{Reduction to atomless laws}
To rule out ties, we first reduce to atomless laws. The event $G=\{x_1+x_2<cx_3\}$ is open in $\R_+^3$. Let $X_1,X_2,X_3$ be iid with law $\mu$, let $U_1,U_2,U_3$ be iid uniform on $[0,1]$ and independent of them, and put $X_i^\varepsilon=X_i+\varepsilon U_i$. The law of $X_i^\varepsilon$ is atomless and lives on $(0,\infty)$. Since $G$ is open, $\mathbf 1_G(X^\varepsilon)\to1$ on $\{X\in G\}$ as $\varepsilon\to0$, so by Fatou's lemma $\liminf_{\varepsilon\to0}\Prob[X^\varepsilon\in G]\ge\Prob[X\in G]$. So it suffices to bound atomless laws on $(0,\infty)$. For those, finitely many sample points are almost surely distinct. Almost surely, no relation $x_i+x_j=cx_l$ between distinct points holds either, because given the other points such a relation pins one point to a single value.

\subsubsection{Triples}
Look at three sample points $u_1<u_2<u_3$. Which of them, multiplied by $c$, can beat the sum of the other two? The smallest never can, since $u_2+u_3>2u_1>cu_1$. Write
\[
  \mathsf T=\{u_1+u_2<cu_3\}\ \text{(the top point wins)},\qquad \mathsf M=\{u_1+u_3<cu_2\}\ \text{(the middle point wins)}.
\]
If the middle point wins, so does the top point: adding the strict inequalities $u_2-u_3<0$ and $c(u_2-u_3)<0$ to $u_1+u_3-cu_2<0$ gives $u_1+u_2-cu_3<0$. So every triple is of one of three kinds: \emph{both} ($\mathsf T$ and $\mathsf M$), \emph{top only}, or \emph{neither}.

\subsubsection{The square}
This is the step that uses independence. Let $X_1,\dots,X_4$ be iid with an atomless law on $(0,\infty)$, and call the first two points $R=\{X_1,X_2\}$ the roots. Adding a point $z\notin R$ gives a triple. Its \emph{configuration} records the kind of the triple and the position of $z$:
\[
 \kappa(R,z)\in\mathcal K
 :=\{\text{both},\text{top only},\text{neither}\}
 \times\{\text{top},\text{middle},\text{bottom}\},
\]
so $|\mathcal K|=9$. A \emph{score} is a function $v:\mathcal K\to\R$, equivalently a vector $v\in\R^{\mathcal K}\cong\R^9$. Let $Z$ be a fresh independent point. For a two-element set $s\subset(0,\infty)$ put $g(s)=\E[v(\kappa(s,Z))]$. Then $g(R)=\E[v(\kappa(R,Z))\mid R]$ is the mean score given the roots. Given $R$, the points $X_3$ and $X_4$ are independent copies of $Z$, so
\begin{equation}\label{eq:square}
  \E\bigl[v(\kappa(R,X_3))\,v(\kappa(R,X_4))\bigr]\;=\;\E\bigl[g(R)^2\bigr]\;\ge\;0 .
\end{equation}
The same identity has a matrix form, which is the one the general method uses. Let $h(s)\in[0,1]^{\mathcal K}$ have coordinates $h_\kappa(s)=\Prob[\kappa(s,Z)=\kappa]$, which sum to $1$. Then $g(s)=v^\top h(s)$, and
\[
 A:=\E[h(R)h(R)^\top]\in\R^{\mathcal K\times\mathcal K}
 \cong\R^{9\times9},\qquad
 \E[g(R)^2]=v^\top A v\ge0.
\]
In particular, $A$ is symmetric positive semidefinite. We use the following score vector, arranged as a $3\times3$ table by kind and position. Its entries look arbitrary at this point; they come from a small semidefinite program, as explained after Proposition~\ref{prop:example}.
\begin{center}
\begin{tabular}{lccc}
\toprule
kind of triple & new point on top & in the middle & at the bottom\\
\midrule
both & $6/11$ & $6/11$ & $7/11$\\
top only & $-2/11$ & $-4/11$ & $-4/11$\\
neither & $-7/11$ & $-7/11$ & $-7/11$\\
\bottomrule
\end{tabular}
\end{center}

\subsubsection{Patterns}
To use \eqref{eq:square}, we compute both the objective and the square from a finite amount of information about the sample, its pattern. Sort the four points as $x_1<x_2<x_3<x_4$. The \emph{pattern} of the sample records, for each of the eight relations
\begin{align*}
  &0<cx_3-x_1-x_2,\quad 0<cx_4-x_1-x_2,\quad 0<cx_2-x_1-x_3,\quad 0<cx_4-x_1-x_3,\\
  &0<cx_2-x_1-x_4,\quad 0<cx_3-x_1-x_4,\quad 0<cx_4-x_2-x_3,\quad 0<cx_3-x_2-x_4,
\end{align*}
the sign of its right-hand side ($+$ when the relation holds, $-$ otherwise). Equalities have probability zero under an atomless law. The remaining four relations $x_i+x_j<cx_l$ with $l<\min(i,j)$ never hold. In the gap coordinates $d_1=x_1$, $d_a=x_a-x_{a-1}$, each relation is a strict homogeneous linear inequality in $d\in(0,\infty)^4$. So a sign vector is realizable if and only if an open polyhedral cone is nonempty. For $c=6/5$ exactly $14$ of the $2^8$ sign vectors are realizable; Table~\ref{tab:example} lists them. A depth-first search proves that no others occur. It abandons each impossible partial sign vector with a Gordan-type certificate: a nonnegative combination of its prescribed strict inequalities that is coordinatewise of the wrong sign (Appendix~\ref{app:gordan}). The implication $\mathsf M\Rightarrow\mathsf T$ above is the simplest example.

Both sides of our inequality are averages over patterns. Write $x_P$ for the probability that the sample has pattern $P$.
\begin{itemize}
  \item The event: by exchangeability, $\Prob[X_1+X_2<cX_3]=\sum_Px_P\,\obj(P)$, where $\obj(P)$ is the fraction of the twelve choices (pair, third point) for which the pair loses to the third point, that is, the number of $+$ signs in $P$ divided by $12$.
  \item The square: by exchangeability we may choose the roots $R$ and the order of the two remaining points uniformly at random among the twelve choices. Then the left side of \eqref{eq:square} is $\sum_Px_P\,m(P)$, where $m(P)$ is the average of $v(\kappa(R,a))\,v(\kappa(R,b))$ over these twelve choices, computed in the pattern $P$. For example, in the pattern \verb|--------| every triple is of kind \emph{neither}, so $m=(7/11)^2=49/121$.
\end{itemize}

\subsubsection{The bound}
Table~\ref{tab:example} gives, for each pattern, the objective, the square term and their sum $\val(P)=\obj(P)+m(P)$.

\begin{table}[ht]
\centering
\small
\begin{tabular}{lccc@{\qquad}lccc}
\toprule
pattern & $12\,\obj$ & $m$ & value & pattern & $12\,\obj$ & $m$ & value\\
\midrule
\verb|++++++--| & 6 & $-13/242$ & $\mathbf{54/121}$ & \verb|++-+-+--| & 4 & $-37/726$ & $205/726$\\
\verb|++++-+--| & 5 & $-35/363$ & $155/484$ & \verb|++-+--+-| & 4 & $34/363$ & $155/363$\\
\verb|++++--+-| & 5 & $5/363$ & $625/1452$ & \verb|++-+----| & 3 & $56/363$ & $587/1452$\\
\verb|++++----| & 4 & $-1/726$ & $241/726$ & \verb|++------| & 2 & $5/22$ & $13/33$\\
\verb|++-+-+++| & 6 & $-43/726$ & $160/363$ & \verb|-+-+--+-| & 3 & $15/121$ & $181/484$\\
\verb|++-+-++-| & 5 & $-8/121$ & $509/1452$ & \verb|-+-+----| & 2 & $149/726$ & $45/121$\\
 & & & & \verb|-+------| & 1 & $217/726$ & $185/484$\\
 & & & & \verb|--------| & 0 & $49/121$ & $49/121$\\
\bottomrule
\end{tabular}
\vspace{6pt}\caption{The 14 realizable patterns of four points for $c=6/5$, in the order of the eight relations listed above, with $\val(P)=\obj(P)+m(P)$. The largest value is $54/121$. Appendix~\ref{app:gordan} lists the certificates showing that no other sign vector occurs.}
\label{tab:example}
\end{table}

\begin{proposition}\label{prop:example}
$S_2(6/5)\le 54/121<0.4463$.
\end{proposition}

\begin{proof}
By \eqref{eq:square}, for every atomless law,
\[
  \Prob[X_1+X_2<\tfrac65X_3]=\sum_Px_P\,\obj(P)\;\le\;\sum_Px_P\,\bigl(\obj(P)+m(P)\bigr)\;\le\;\max_P\val(P)=\frac{54}{121},
\]
using $\sum_Px_P=1$ and Table~\ref{tab:example}.
\end{proof}

Adding $m(P)$ is what helps. Without it, the best we can say is $\max_P\obj(P)=1/2$, the four-point bound of \citet{BF}. Both patterns in which six of the twelve choices succeed have $m(P)<0$, and $m(P)$ is largest on the patterns with the fewest successes. Since the square has nonnegative mean, no law can put all its weight on the good patterns. It must balance them with patterns where $m$ is positive and the objective is small.

Each entry of Table~\ref{tab:example} can be checked by hand, and so can the completeness of the list (Appendix~\ref{app:gordan}). So Proposition~\ref{prop:example} can be verified with pencil and paper. Our verifier also checks everything in exact arithmetic. We found $v$ by solving a small semidefinite program. For these flags its optimal matrix has rank one, so the optimal certificate is itself a single square, and $v$ is its vector rounded to elevenths. The optimum of that program is about $0.43568$. With six points the certified bound is $0.4305620$, and recording a pair color gives the exact value $(2\sqrt7-4)/3$ with three points (Section~\ref{sec:plateau}).

\subsection{The general method}\label{sec:general}
The general method repeats the example with more bookkeeping: a finite list of patterns, the objective as an average over patterns, squares with nonnegative mean, and the bound by the worst pattern. We consider
\[
  \sup_\mu \Prob_{\mu^{\otimes n}}\bigl[(X_1,\dots,X_n)\in G\bigr]
\]
where, for an integer $m\ge1$ and fixed coefficient vectors $a_1,\dots,a_m\in\R^n$,
\[
 G=\{x\in\R_+^n:a_j^\top x<0\text{ for }j=1,\dots,m\}.
\]
For the four-roll game, $n=4$ and $m=1$, with the single vector $a_1=(1,1,1,-2)^\top$. So $a_1^\top x<0$ is exactly $x_1+x_2+x_3<2x_4$. Fix an integer sample size $s\ge n$. For a fixed law $\mu$, we abbreviate $\Prob_{\mu^{\otimes n}}[(X_1,\dots,X_n)\in G]$ to $\Prob[G]$, where $X_1,\dots,X_n$ are iid with law $\mu$.
\begin{enumerate}
  \item \emph{Patterns:} The pattern of a sorted sample $x_1<\dots<x_s$ is the sign vector of the dot products $a_j^\top(x_{i_1},\dots,x_{i_n})$, for $j=1,\dots,m$ and all ordered $n$-tuples of distinct sample indices. Forms of constant sign on sorted samples are dropped. Gordan certificates prove that a finite list $\cP_s$ contains every realizable pattern. The list may also contain unrealizable patterns. This keeps the bound valid but can weaken it, since the semidefinite constraints need not detect that these patterns cannot occur. A pattern may be removed once its non-realizability has been certified. For an atomless law $\mu$, let $x_P$ be the probability of pattern $P$. Then $x=(x_P)_{P\in\cP_s}\in[0,1]^{\cP_s}$ and $\sum_Px_P=1$. This probability vector depends on $\mu$.
  \item \emph{Objective:} $\Prob[G]=\sum_Px_P\obj(P)$, where the fixed function $\obj:\cP_s\to[0,1]$ assigns to $P$ the fraction of ordered $n$-tuples of distinct sample points that lie in $G$.
  \item \emph{Types and flags:} For integers $0\le t<f$ with $2f-t\le s$, a \emph{type} $\sigma$ is a pattern on $t$ points, and a \emph{$\sigma$-flag} is a pattern on $f$ points together with $t$ of its positions, the \emph{roots}, on which it restricts to $\sigma$. Formally, with $[f]=\{1,\dots,f\}$,
  \begin{align*}
    \text{type: }&\sigma\in\cP_t,\\
    \text{$\sigma$-flag: }&F\in\mathcal F_\sigma
    :=\{(P,R):P\in\cP_f,\ R\subseteq[f],\ |R|=t,\ P|_R=\sigma\}.
  \end{align*}
  Because the points are ordered by value, their ranks identify them uniquely, and the roots inherit their labels from this order. A set of $f$ sample points, $t$ of them marked as roots, \emph{carries} the flag $F$ if its pattern, with the marked points as roots, is $F$. In the example, $t=2$ and $f=3$, the roots are the pair $R$, and the nine configurations $\kappa$ are the flags of the unique type.

  For each type $\sigma$ the fixed data are symmetric matrices $M_\sigma(P)\in\R^{\mathcal F_\sigma\times\mathcal F_\sigma}$, one for each pattern $P$, with rows and columns indexed by flags:
  \[
    M_\sigma(P)[F_1,F_2]=\frac{\#\{(R,A,B):\ P|_R=\sigma,\ R\cup A\text{ carries }F_1,\ R\cup B\text{ carries }F_2\}}{\binom st\binom{s-t}{f-t}\binom{s-f}{f-t}}.
  \]
  The count runs over pairwise disjoint subsets $R,A,B$ of the $s$ points with $|R|=t$ and $|A|=|B|=f-t$; the denominator is the number of all such triples, whatever their type, and the roots of the flags are the points of $R$. An individual $M_\sigma(P)$ need not be positive semidefinite; Lemma~\ref{lem:psd} shows that its expectation $\sum_Px_PM_\sigma(P)$ under an iid sample is. The variables are symmetric matrices $Q_\sigma$ of the same size, constrained to be positive semidefinite, and the pairing is $\langle Q,M\rangle:=\operatorname{tr}(Q^\top M)=\sum_{F,F'}Q_{F,F'}M_{F,F'}$. For every such $Q_\sigma$ we will have $\sum_Px_P\,\langle Q_\sigma,M_\sigma(P)\rangle\ge0$ (Lemma~\ref{lem:psd}).
  \item \emph{Certificate:} Then $\sup\Prob[G]\le\lambda=\max_{P\in\cP_s}\bigl(\obj(P)+\sum_\sigma\langle Q_\sigma,M_\sigma(P)\rangle\bigr)$.
\end{enumerate}

In probabilistic terms, a type is the pattern of the conditioning points, a flag is the pattern of one extra batch together with where the conditioning points sit, and the next lemma says that, given the conditioning points, two disjoint batches are iid.

\begin{lemma}\label{lem:psd}
Let $\mu$ be atomless, let $Y_R$ be $t$ iid points with law $\mu$, and let $h(Y_R)\in[0,1]^{\mathcal F_\sigma}$ be the column vector with entries $h_F(Y_R)=\Prob[\,Y_R\cup Y_A\text{ carries }F\mid Y_R\,]$, where $Y_A$ is an independent batch of $f-t$ further points and $F\in\mathcal F_\sigma$. When the roots have type $\sigma$, the coordinates sum to $1$; otherwise they are all zero. Then
\[
  \sum_Px_P\,M_\sigma(P)\;=\;\E\bigl[\mathbf 1\{Y_R\text{ has type }\sigma\}\,h(Y_R)h(Y_R)^\top\bigr],
\]
and hence $\sum_Px_P\langle Q_\sigma,M_\sigma(P)\rangle=\E\bigl[\mathbf 1\{Y_R\text{ has type }\sigma\}\,h(Y_R)^\top Q_\sigma h(Y_R)\bigr]\ge0$ for every positive semidefinite $Q_\sigma$.
\end{lemma}

\begin{proof}
Take $s$ iid points and, independently, a uniformly random triple $(R,A,B)$ of disjoint sets of labels as in the definition of $M_\sigma$. The points indexed by $R$, $A$ and $B$ are then three independent batches, distributed as $Y_R$, $Y_A$ and an independent copy $Y_B$ of $Y_A$. Conditionally on the values, the positions that $R$, $A$, $B$ occupy in the sorted sample form a uniformly random triple, because the labels were chosen uniformly and independently of the values. Whether $R$ has type $\sigma$ and which flags $R\cup A$ and $R\cup B$ carry depends only on the pattern of the $s$ points and on these positions, since the pattern of a subset is the restriction of the pattern. Averaging over patterns gives the left side. On the other hand, given $Y_R$ the flags of $Y_R\cup Y_A$ and $Y_R\cup Y_B$ are independent with the same conditional law, which gives the right side.
\end{proof}

The matrix identity averages outer products of conditional probability vectors, just as for the nine-coordinate vector in Section~\ref{sec:example}. In the language of \cite{Razborov07}, the realizable patterns are the models of a universal theory with a linear order and one predicate for each linear form, $M_\sigma(P)$ is the coefficient of $P$ in the level-$s$ expansion of $\llbracket F_1F_2\rrbracket_\sigma$, and Lemma~\ref{lem:psd} is the inequality $\llbracket f^2\rrbracket_\sigma\ge0$.

The pattern matrices $M_\sigma(P)$ and objective values $\obj(P)$ are fixed data. The semidefinite program chooses the matrices $Q_\sigma$ and a scalar $\lambda\in\R$ to
\[
 \text{minimize }\lambda\quad\text{subject to}\quad
 Q_\sigma\succeq0,\qquad
 \obj(P)+\sum_\sigma\langle Q_\sigma,M_\sigma(P)\rangle\le\lambda
 \quad(P\in\cP_s).
\]
Here $Q_\sigma\succeq0$ means positive semidefinite. The pattern probabilities $x_P$ are not variables of this program. We solve it numerically, round $Q_\sigma=L_\sigma L_\sigma^\top$ with rational $L_\sigma$, so that exact positive semidefiniteness is automatic, as is standard \cite{PeyrlParrilo,FRV13}, and compute $\lambda$ exactly. In Proposition~\ref{prop:example}, $s=4$, $t=2$ and $f=3$. There is a unique type $\sigma$, since no comparison involves fewer than three distinct points. The certificate uses $Q_\sigma=vv^\top$, where $v$ is the rounded score vector from the example, so that $m(P)=v^\top M_\sigma(P)v=\langle vv^\top,M_\sigma(P)\rangle$. This matrix is positive semidefinite and rank one; we do not claim it is optimal.

Optimality is not needed for the proof. Any positive semidefinite matrices $Q_\sigma$ give an upper bound for every law $\mu$, namely the maximum certificate value over the complete pattern list $\cP_s$. Optimization only makes the bound smaller. Validity requires three checks: positive semidefiniteness, completeness of the pattern list, and the claimed bound on every listed pattern.

%======================================================================
\subsection{Relation to upper bounds by exchangeability}\label{sec:relation}
%======================================================================

Without the squares, the method gives exactly the bound of \citet{BF}. The squares can only improve it, and both versions converge to the true value as the sample grows. For a sample size $s$, let $\mathrm{MFS}_s$ (for maximum feasible subsystem) be the largest fraction of instances of the inequality that $s$ actual points can satisfy simultaneously, that is, the maximum of $\obj(P)$ over \emph{realizable} patterns. This is the upper bound of \cite[Lemma~2.3]{BF} with $m=s$. Let $\lambda_s$ be the best bound of Section~\ref{sec:general} with sample size $s$, all admissible types and flags, and $\cP_s$ the set of realizable patterns. (Some of our computational certificates use an \emph{outer cover}: a certified list that contains every realizable pattern and possibly some unrealizable ones. This keeps the bound valid and can only weaken it.)

\begin{proposition}\label{prop:sandwich}
For an event $G\subseteq\R_+^n$ as in Section~\ref{sec:general} and every $s\ge n$,
\[
  \sup_\mu\Prob_{\mu^{\otimes n}}[G]\;\le\;\lambda_s\;\le\;\mathrm{MFS}_s\;\le\;\frac{s^n}{s(s-1)\cdots(s-n+1)}\,\sup_\mu\Prob_{\mu^{\otimes n}}[G].
\]
The bound $\lambda_s$ is also nonincreasing in $s$, and $\lambda_s$ converges to the true value as $s\to\infty$.
\end{proposition}

\begin{proof}
The first inequality is Section~\ref{sec:general}. The second holds because $Q=0$ is admissible. The third is \cite[Theorem~2.4 and Remark~2.5]{BF}: the uniform law on an optimal configuration of $s$ points achieves at least $\mathrm{MFS}_s$ times the probability $s(s-1)\cdots(s-n+1)/s^n$ that $n$ independent draws from it are distinct, which is Freedman's bound for sampling with and without replacement \cite{Freedman77}; this factor is at least $(1-n/s)^n\to1$. Monotonicity in $s$: a certificate at level $s$ is one at level $s+1$, because both $\obj$ and $M_\sigma$ at level $s+1$ are averages of their values on the $s$-point restrictions.
\end{proof}

This is the analogue of the familiar fact that the largest induced density of a fixed graph over $s$-vertex graphs decreases to its limit, with the same blow-up argument. The guaranteed rate is far from what happens in practice: for $C$ it only gives $\lambda_7\le\frac{2401}{840}C<2.86\,C$. Table~\ref{tab:levels} shows the actual effect of the semidefinite constraints.

\begin{table}[ht]
\centering
\begin{tabular}{lllll}
\toprule
problem & sample size & $\mathrm{MFS}$ (no PSD) & flag bound & best construction\\
\midrule
$C$ & $s=6$ & $7/15=0.4667$ & $0.4025452$ & $0.4006953$\\
$C$ & $s=7$ & $16/35=0.4571$ & $0.4007545$ & $0.4006953$\\
$C$ & $m=20$ \cite{BF} & $0.4167$ & & \\
$C$ & $m=25$ \cite{BF} (no certificate) & $0.4130$ & & \\
$S_2(6/5)$ & $s=4$ & $1/2$ & $0.4356824^{*}$ & $0.4305008$\\
$S_2(6/5)$ & $s=6$ & $29/60=0.4833$ & $0.4305620$ & $0.4305008$\\
\bottomrule
\end{tabular}
\vspace{6pt}\caption{The bound of \cite{BF} ($\mathrm{MFS}$) and the flag bound at the same sample size; the constructions give $C\ge C_{\mathrm{BF}}\approx0.4006953723$. ${}^*$Numerical optimum; the rounded rank-one certificate of Proposition~\ref{prop:example} gives $54/121\approx0.446$.}
\label{tab:levels}
\end{table}

\begin{remark}[Exchangeability and limit objects]\label{rem:exch}
The pattern law $(x_P)$ of $(X_1,\dots,X_s)$ is defined for any exchangeable sequence without ties, and the identity $\Prob[G]=\sum_Px_P\obj(P)$ uses only exchangeability. The extreme exchangeable laws are urns, uniform orderings of $s$ fixed points, and the pattern law of an urn is a point mass. So the level without semidefinite constraints, $\mathrm{MFS}_s$, is exactly the supremum over finitely exchangeable sequences. This distinction between exchangeable laws and mixtures of iid laws is the starting point of Lacker's argument \cite{Lacker}; see also \cite{DiaconisFreedman80}. Lemma~\ref{lem:psd} uses more: under an iid law, two disjoint batches are independent given the roots. Under a mixture of iid laws they are independent only given the roots and the directing law, but the inequalities of Lemma~\ref{lem:psd} still hold, by averaging over the directing law. By de Finetti's theorem these mixtures are exactly the exchangeable sequences that extend to infinite ones, and since $\Prob[G]$ is linear in the law, mixtures do not raise the supremum. For urns the inequalities can fail. The urn on the four points of the pattern \verb|++++++--| of Table~\ref{tab:example} has $m=-13/242<0$, so it violates \eqref{eq:square}.

Our relaxation keeps only finitely many of these constraints, so the pattern laws it allows need not come from any law $\mu$. Even the limit objects of the theory, the consistent pattern laws satisfying all of the constraints, need not (compare \cite[Proposition~3]{Goaoc}). For our objectives this makes no difference in the limit. A limit object is a distribution on realizable $s$-point patterns for every $s$, so its value is at most $\mathrm{MFS}_s$, and by Proposition~\ref{prop:sandwich} the supremum over limit objects equals the supremum over laws. The same principle, that independent copies of a source are exchangeable and conditionally independent, drives the inflation technique \cite{WSF19}.
\end{remark}

%======================================================================
\subsection{The four-roll game with six points: \texorpdfstring{$C\le0.4025452$}{C <= 0.4025452}}\label{sec:bta}
%======================================================================

We now run the method on the four-roll game, recording additive relations only. Section~\ref{sec:four-exact} supersedes these bounds with the exact value. We keep them because they show how far the plain method gets, and they fill the first rows of Table~\ref{tab:levels}.

Here $n=4$ and the instances are $x_i+x_j+x_k<2x_l$. On a sorted sample an instance can hold only if $l$ is larger than the middle index of $\{i,j,k\}$, since otherwise $x_j+x_k\ge2x_l$. So every four points $a<b<c<d$ give two candidate instances, $x_a+x_b+x_c<2x_d$ and $x_a+x_b+x_d<2x_c$, and the second implies the first. That makes $2\binom s4$ candidate instances.

Six points give $30$ candidate instances and $1{,}290$ realizable patterns. A depth-first search over sign vectors proves that this list is complete, and each of its $6{,}961$ pruned branches carries an exact Gordan certificate. No instance involves fewer than four points, so there is a unique type on two points and a unique type on three, and we use each with flags on four points. We also use the three types on four points, with flags on five points. The blocks have $18$, $12$, $51$, $63$ and $36$ flags. The semidefinite program has value $0.4025451$, and the rounded certificate gives
\begin{equation}\label{eq:s6}
  C\;\le\;\frac{201272564718335614989}{500000000000000000000}\;<\;0.4025452 .
\end{equation}

\subsection{The four-roll game with seven points: \texorpdfstring{$C\le0.4007545$}{C <= 0.4007545}}\label{sec:seven}
Seven points give $70$ candidate instances and $136{,}800$ realizable patterns. At this size we wanted two independent proofs that the list is complete.
\begin{enumerate}
  \item Directly, by the depth-first search, with $1{,}557{,}094$ Gordan certificates. This number depends on the order in which the relations are decided; ours is recorded with the certificates.
  \item By extension. Every six-point restriction of a realizable seven-point pattern is realizable, so it lies in the certified six-point list. Joining the seven restrictions gives $344{,}093$ candidates. Gordan certificates exclude $207{,}293$ of them, and $136{,}800$ remain.
\end{enumerate}
The two lists coincide.

The certificate uses the families $(t,f)=(1,4),(3,5),(5,6)$. The first two have $12$ and $300$ flags; the third has thirty types with $7{,}740$ flags in all. This is the full level-seven relaxation. Families with $2f-t<7$ add nothing: a flag of such a family can be enlarged by an extra root (or extra non-root points) that the score ignores, which turns its constraint into one of these families.

With $136{,}800$ patterns the program is large, so we solved it by column generation. An interior-point method runs on a working set of patterns, and each round adds the most violated patterns. We stopped when the certificate, evaluated over all $136{,}800$ patterns, was within $5\cdot10^{-8}$ of the working-set optimum. The final working set has $4{,}454$ patterns. The matrices are given as $Q_\sigma=L_\sigma L_\sigma^\top/4^{40}$ with integer $L_\sigma$, and the bound
\begin{equation}\label{eq:s7}
  C\;\le\;\frac{33913767225127410008898759}{84624807373024042229432320}\;<\;0.4007545
\end{equation}
is computed over all $136{,}800$ patterns in integer arithmetic over a common denominator. The numerical optimum is about $0.400754436$.

Going from six to seven points shrinks the gap above $C_{\mathrm{BF}}$ from about $1.8\cdot10^{-3}$ to less than $5.92\cdot10^{-5}$, but does not close it. Tomar's nine-point certificate gives $C\le0.40069634$ \cite{Tomar}. Section~\ref{sec:four-exact} closes the gap in a different way: it stays at four points and records more relations.

\subsection{Verification}
Each bound is checked in exact arithmetic by a verifier that shares no code with the search. For six points the verifier is about $100$ lines of Python. It replays all completeness certificates, recomputes every $M_\sigma(P)$ from the definition, and outputs \eqref{eq:s6} in a few seconds. For seven points, a pure-Python check replays all $1{,}557{,}094$ completeness certificates and confirms that the surviving sign vectors are exactly the $136{,}800$ patterns on which a stand-alone program evaluates \eqref{eq:s7}. This takes a few minutes. A second verifier, written independently of the search, re-derives the pattern list by the extension route and also confirms \eqref{eq:s7}. The patterns, flags and matrices were also reimplemented independently, with agreement pattern by pattern. As sanity checks, the verifiers return $7/15$ and $16/35$, the $\mathrm{MFS}$ values of Table~\ref{tab:levels}, when all $Q_\sigma$ are set to zero, and the six-point bound agrees with the numerical optimum to seven digits.

\begin{remark}[Rationality]
A certificate with rational entries gives a rational bound, whereas $C_{\mathrm{BF}}$ is irrational\footnote{Clearing denominators in \eqref{eq:lower}, $q=C_{\mathrm{BF}}$ is the value at which the cubic $qp^3+(q+1)p^2+(q-2)p+q=(1+p+p^2+p^3)\bigl(q-p(2-p)/(1+p+p^2+p^3)\bigr)$ acquires a double root in $(0,1)$. Setting its discriminant to zero gives $16q^4+8q^3+99q^2-32q-4=0$, which has no rational roots.}. So no single rational certificate can prove the exact value. Algebraic entries remove this obstruction. The colored four-point certificate of Section~\ref{sec:four-exact} is sharp over $\mathbb Q(r)$, where $r\in(0,1)$ is the root of the quartic defined there. For the additional exact $S_2$ plateaus, Section~\ref{sec:thmD} combines algebraic certificates with a fixed-point inequality.
\end{remark}

%======================================================================
\section{Improved flag SDP by recording pair colors}\label{sec:pair-colors}
%======================================================================
The relations recorded in a pattern need not be limited to those defining
the objective. Pairwise comparisons can expose information that the
additive relations alone do not record on a small sample. We first give
an exact three-point certificate for the first plateau of $S_2$, then
describe the refinement for the general method of Section~\ref{sec:general}.

\subsection[Pair colors on a two-against-one example]{Pair colors on a two-against-one example: Exact value of \texorpdfstring{$S_2(c)$ for $c\in(1,\sqrt2\,]$}{S2(c) for 1<c<=sqrt(2)}}\label{sec:plateau}

\begin{theorem}\label{thm:plateau}
For every $1<c\le\sqrt2$,
\[
 S_2(c)=V:=\frac{2\sqrt7-4}{3}=0.4305008\ldots.
\]
The upper bound has an optimal three-point flag certificate with one
pair-color cutoff and a single rank-one positive semidefinite matrix.
\end{theorem}

\begin{proof}
By the atomless approximation in Section~\ref{sec:example}, it suffices
for the upper bound to consider atomless laws on $(0,\infty)$.
Choose $r$ with $c\le r\le2/c$, and color a pair $x<y$ according to
whether $y<rx$; call such a pair an edge. The resulting ordered graph
has the interval property: an edge $xz$, for $x<y<z$, implies both
$xy$ and $yz$. There are therefore five three-point graph patterns:
$\mathsf0$ (no edges), $\mathsf L$ (only the lower pair), $\mathsf H$
(only the upper pair), $\mathsf P$ (both consecutive pairs), and
$\mathsf T$ (a triangle).

For a sorted triple $x<y<z$, let $P$ denote its ordered edge/nonedge
pattern, one of $\mathsf0,\mathsf L,\mathsf H,\mathsf P,\mathsf T$,
and let $N$ be the number of winning target positions:
\[
 N=\mathbf1\{x+y<cz\}+\mathbf1\{x+z<cy\}.
\]
A middle win implies $z<cy\le ry$ and
$x<(c-1)y\le y/r$, since $r(c-1)\le2(c-1)/c\le1$.
It can therefore occur only in $\mathsf H$.
A triangle gives $x>z/r$, so $x+y>2z/r\ge cz$ and there are no wins.
Consequently the additive objective satisfies
\begin{equation}\label{eq:pair-majorant}
 \frac{N}{3}\le g(P):=
 \frac{1+\mathbf1_{\mathsf H}(P)-\mathbf1_{\mathsf T}(P)}3.
\end{equation}

Here the type consists of a single marked point, the root. There is
only one such type: one point carries neither a pair color nor an additive
comparison. A flag adds one point to this root, recording both the color
of their pair and which of the two positions is the root. Thus the family
$(t,f)=(1,2)$ has exactly four flags. Writing $x$ for the root and $y$ for
the additional point, they are
\[
\begin{array}{c|c|c|l}
\text{flag}&\text{ordered pair (root boxed)}&\text{pair color}&\text{condition}\\\hline
F_1&y<\boxed{x}&\text{nonedge}&0<y<x/r\\
F_2&y<\boxed{x}&\text{edge}&x/r<y<x\\
F_3&\boxed{x}<y&\text{edge}&x<y<rx\\
F_4&\boxed{x}<y&\text{nonedge}&rx<y.
\end{array}
\]
The boundary equalities have probability zero for atomless laws. The five graphs
$\mathsf0,\mathsf L,\mathsf H,\mathsf P,\mathsf T$ are the three-point
patterns on which the certificate is checked; the rows and columns of
its single PSD matrix are indexed by $F_1,F_2,F_3,F_4$.
For a fixed root $x$, the conditional flag probabilities are
\[
 h(x)=\bigl(\mu((0,x/r)),\ \mu((x/r,x)),\
             \mu((x,rx)),\ \mu((rx,\infty))\bigr)^\top.
\]
Two independent additional points give the matrix $h(x)h(x)^\top$.
Averaging over the root yields $\E M(P)=\E[h(X)h(X)^\top]\succeq0$,
which is Lemma~\ref{lem:psd} for this particular flag family.

Concretely, $M_{ij}(P)$ is $1/6$ times the number of ordered choices
$(u,a,b)$ of the three distinct vertices for which $u$ is the root,
$(u,a)$ carries $F_i$, and $(u,b)$ carries $F_j$.

The SDP minimizes a scalar upper bound $\lambda$ over one symmetric
matrix $Q\in\mathbb R^{4\times4}$. The cutoff $r$, the five patterns,
and their matrices $M(P)$ are fixed. Each pattern contributes the linear
constraint $g(P)+\langle Q,M(P)\rangle\le\lambda$; explicitly, the program is
\begin{equation}\label{eq:pair-plateau-sdp}
\begin{aligned}
 &\operatorname{minimize}_{\lambda\in\mathbb R,\ Q=Q^\top}\quad\lambda\\
 &\text{subject to}\quad Q\succeq0,\\
 &\qquad 1+Q_{11}+Q_{14}+Q_{44}\le3\lambda
 &&(\mathsf0),\\
 &\qquad 1+Q_{11}+Q_{24}+Q_{34}\le3\lambda
 &&(\mathsf L),\\
 &\qquad 2+Q_{12}+Q_{13}+Q_{44}\le3\lambda
 &&(\mathsf H),\\
 &\qquad 1+Q_{12}+Q_{23}+Q_{34}\le3\lambda
 &&(\mathsf P),\\
 &\qquad Q_{22}+Q_{23}+Q_{33}\le3\lambda
 &&(\mathsf T).
\end{aligned}
\end{equation}
Every feasible $\lambda$ bounds $\E[g(P)]$, and therefore $S_2(c)$ by
\eqref{eq:pair-majorant}. The matrix $Q$ is allowed to have any rank.
We now exhibit a rank-one solution with $\lambda=V$; the matching lower
construction below proves that this solution minimizes
\eqref{eq:pair-plateau-sdp}.
Put $\theta=(\sqrt7-1)/2$ and $v=(-\theta,1,1,0)^\top$, and take
\begin{equation}\label{eq:pair-plateau-matrix}
 Q=Vvv^\top=
 \begin{pmatrix}
 2\sqrt7-5&\sqrt7-3&\sqrt7-3&0\\
 \sqrt7-3&V&V&0\\
 \sqrt7-3&V&V&0\\
 0&0&0&0
 \end{pmatrix}\succeq0.
\end{equation}
For example, in $\mathsf H$, with vertices $x<y<z$ and only the edge $yz$,
choosing $x,y,z$ as the root gives the flag pairs
$(F_4,F_4)$, $(F_1,F_3)$, and $(F_1,F_2)$, respectively, before swapping
the two extensions. Thus
\[
 v^\top M(\mathsf H)v
 =\frac{v_4^2+v_1v_3+v_1v_2}{3}=-\frac{2\theta}{3}.
\]
The same count for the other four patterns gives the entire certificate:
\[
\begin{array}{c|ccccc}
P&\mathsf0&\mathsf L&\mathsf H&\mathsf P&\mathsf T\\\hline
v^\top M(P)v&\theta^2/3&\theta^2/3&-2\theta/3&(1-\theta)/3&1\\
g(P)+\langle Q,M(P)\rangle&V&V&V&5V/6&V.
\end{array}
\]
The identities follow from $2\theta^2+2\theta=3$ and
$V=2/(3+2\theta)$. Lemma~\ref{lem:psd} gives
$\E\langle Q,M(P)\rangle\ge0$, and hence
\[
 \Prob[X_1+X_2<cX_3]=\E[N]/3\le\E[g(P)]\le V.
\]
If the additive labels are retained, they merely split these graph
patterns into finer patterns: \eqref{eq:pair-majorant} shows that the
same matrix is a certificate for their actual additive objective.

For the reverse inequality, put masses $q(1-q)^k$ at
$\varepsilon^k$, $k\ge0$, where
\[
 q=\frac{4-\sqrt7}{3},\qquad
 0<\varepsilon<\min\{c/2,c-1,1/c\}.
\]
A largest-scale target wins against two smaller draws; either of two
largest-scale targets wins when the third draw is smaller; three draws
on one scale never win. The success probability $p$ therefore satisfies
\[
 p=(1-q)^3p+q(1-q)^2+2q^2(1-q),
\]
whose solution is $p=V$. This proves the lower bound and also the
optimality of the colored three-point SDP certificate.
\end{proof}

The direction $v$ can be read from the construction. At graph level it
consists of ordered cliques of masses $q(1-q)^k$. A root in the $k$th
clique has neighbor mass $q(1-q)^k$ and lower-nonneighbor mass
$(1-q)^{k+1}$; their ratio is $q/(1-q)=\theta$. Thus the flag combination
$F_2+F_3-\theta F_1$ vanishes on the extremal law.
Equating the certificate values on $\mathsf0,\mathsf H,\mathsf T$
then determines $\theta$ and the multiplier $V$.
All entries of \eqref{eq:pair-plateau-matrix} lie in $\mathbb Q(\sqrt7)$,
and its factorization and the five displayed values give an analytic
proof of exact SDP optimality.

\subsection{Recording pair colors: The general case}\label{sec:pair-colors-general}

Return to the event $G$ and sample size $s$ of Section~\ref{sec:general}.
Fix cutoffs $1<q_1<\cdots<q_m$. In addition to the original relations,
record on each ordered pair $x_i<x_j$ all comparisons
\begin{equation}\label{eq:pair-color-relations}
 x_j<q_\ell x_i,\qquad 1\le\ell\le m.
\end{equation}
Equivalently, the pair receives one of $m+1$ colors according to the
interval containing $x_j/x_i$. We again work first with positive atomless
laws, so the boundaries $x_j=q_\ell x_i$ have probability zero.
The original relations still determine the objective; the colors supply
additional information for the types and flags.

A colored pattern restricts to a colored pattern on every subset.
We therefore define colored types, flags and product matrices
$\widetilde M_\rho$ by exactly the same counting rule as before.
Given the root values, the two disjoint extensions remain independent
with the same conditional law, so Lemma~\ref{lem:psd} applies unchanged.
The resulting program is
\[
 \min\lambda\quad\text{subject to}\quad
 \widetilde Q_\rho\succeq0,\qquad
 \obj(\widetilde P)+
 \sum_\rho\langle\widetilde Q_\rho,
                 \widetilde M_\rho(\widetilde P)\rangle
 \le\lambda
 \quad(\widetilde P\in\widetilde{\cP}_s).
\]
There are more patterns, but also more flag coordinates and more root
types. In particular, additive comparisons between three distinct points
give only one two-point type, whereas the colors allow a separate PSD
matrix for each ratio interval of the two roots. In the plateau example,
they already distinguish four two-point flags at a single root.

\begin{proposition}[Refining the recorded relations]\label{prop:pair-refinement}
Fix a sample size and flag families, and use the realizable pattern sets.
If every flag in these families is refined by the pair colors, the optimal
colored SDP bound is at most the optimal bound without those colors.
\end{proposition}
\begin{proof}
Let $\pi$ forget the colors. For a refined root type $\rho$ lying above
$\sigma=\pi(\rho)$, let $A_\rho$ map refined flag coordinates to coarse
ones by summing all coordinates whose flags have the same image under
$\pi$. Given any coarse certificate $Q_\sigma$, set
\[
 \widetilde Q_\rho=A_\rho^\top Q_\sigma A_\rho\succeq0.
\]
Partitioning the root and extension choices by their refined types gives,
for every colored pattern $\widetilde P$,
\[
 \sum_\rho\langle\widetilde Q_\rho,
                    \widetilde M_\rho(\widetilde P)\rangle
 =\sum_\sigma\langle Q_\sigma,M_\sigma(\pi(\widetilde P))\rangle.
\]
The objective is unchanged by forgetting colors, so this lifts the coarse
certificate with the same upper bound. Taking infima proves the claim.
\end{proof}
The same lifting argument applies to certified outer covers whenever
forgetting colors maps the refined cover into the coarse cover.
Thus increasing the sample size and refining the recorded relations are
two ways to strengthen the relaxation. At a fixed sample size, adding
cutoffs preserves every certificate already available with fewer colors.

\paragraph{Choice of cutoffs and exact verification.}
Candidate constructions suggest ratios at which to split the types.
For the ladders of Section~\ref{sec:two}, useful cutoffs come from
\[
 r_0=1,\quad r_{j+1}=(1+r_j)/c,\qquad
 b_0=4/c^2,\quad b_{j+1}=(1+b_j)/c,
\]
and ratios between these locations, such as $r_2/r_1=(c+2)/(2c)$.
The choice of cutoffs does not restrict the laws under consideration:
every positive configuration induces a coloring, and the upper certificate
is checked on a cover containing all realizable colored patterns.
More cutoffs can improve the bound, at the cost of larger pattern lists
and PSD blocks; Proposition~\ref{prop:pair-refinement} does not compare
programs with different sample sizes.

The comparisons \eqref{eq:pair-color-relations} are homogeneous linear
inequalities, just like the original additive comparisons. Completeness
is therefore certified by the same Gordan argument in positive gap
coordinates, using exact algebraic arithmetic when the cutoffs are
algebraic. The final certificate is checked by evaluating its score on
every retained pattern and verifying positive semidefiniteness.
The first-plateau certificate is especially small: a positive algebraic
scalar times one outer product, with only five graph cases to check.

%======================================================================
\subsection[The exact value of the four-roll supremum]{The exact value of \texorpdfstring{$\sup_\mu\Prob[X_1+X_2+X_3<2X_4]$}{sup mu P[X1+X2+X3 < 2X4]}}\label{sec:four-exact}
%======================================================================
Three pairwise ratio cutoffs give an exact certificate on four points.
Let $r$ be the unique root in $(0,1)$ of
$r^4-4r^3-3r^2-2r+2=0$, and put
\[
 B=C_{\mathrm{BF}}=\frac{r(2-r)}{1+r+r^2+r^3}
   =\frac{r^3-6r^2+6r}{4}.
\]
The last identity follows from the equation for $r$.

For four sorted positive points $x_1<x_2<x_3<x_4$, record the two
possible successes
\[
 x_1+x_2+x_3<2x_4,\qquad x_1+x_2+x_4<2x_3,
\]
and, for each pair $i<j$, the three comparisons
$x_j<q x_i$, with $q\in\{7/6,3/2,2\}$.
These twenty signs determine the pattern $P$ and its objective
$\obj(P)$, the number of successes divided by four.
Boundary equalities have probability zero under an atomless law.

Use the flag family $(t,f)=(2,3)$: two roots and one additional point.
There are four root types, according to the interval containing the
ratio of the larger root to the smaller:
$(1,7/6)$, $(7/6,3/2)$, $(3/2,2)$, and $(2,\infty)$.
A flag records the pair comparisons on three points and the positions
of its two roots.
For each type $\sigma$, let $M_\sigma(P)$ be the matrix defined in
Section~\ref{sec:general}; here its denominator is
$\binom42\binom21\binom11=12$.

The exact certificate in \path{anc/proofs/s3/0.400695372__c2/cert_exact.json} gives
four matrices $Q_\sigma$ over $\mathbb Q(r)$ satisfying
\begin{equation}\label{eq:direct-four-cert}
 Q_\sigma\succeq0,\qquad
 \obj(P)+\sum_\sigma\langle Q_\sigma,M_\sigma(P)\rangle\le B
 \quad(P\in\cP_4).
\end{equation}
The cover $\cP_4$ has $259$ retained patterns; $871$ excluded prefixes,
each with a rational Gordan witness, certify its completeness.
The matrices are stored as $Q_\sigma=W_\sigma Y_\sigma W_\sigma^\top$,
with $Y_\sigma$ positive definite. Their sizes are
\[
\begin{array}{c|rrrr}
\text{root ratio interval}&(1,7/6)&(7/6,3/2)&(3/2,2)&(2,\infty)\\ \hline
\text{stored flag coordinates}&15&20&19&21\\
\dim Y_\sigma&9&13&14&8
\end{array}
\]
Unlisted flag coordinates have zero rows and columns.
Exactly $57$ retained patterns have equality in
\eqref{eq:direct-four-cert}; every other one has slack greater than
$0.0223$.

Why these three cutoffs? They come from the lower-bound construction of \citet{BF}. On one scale it places the atoms $1-2^{-k}$, $k=1,\dots,N$. The atom at $0$ has mass $r$, the same $r$ as in $B$, and is replaced by a copy of the whole law at a much smaller scale. Two points on the same scale have ratio less than $2$, while two points on different scales have a much larger ratio. So the comparison with $2$ records whether two points share a scale. This is how a four-point certificate reflects the separation of scales. The cutoffs $3/2$ and $7/6$ are the ratios $\frac34:\frac12$ and $\frac78:\frac34$ of the first atoms on a scale. They split pairs on one scale by how close the two points are.

A sharp certificate must also fit the construction: near-optimal laws have to concentrate on the patterns without slack (Remark~\ref{sec:structure}). We checked this numerically. In a simulation of the construction with $58$ atoms per scale, $99.2\%$ of the samples have one of the $57$ tight patterns ($36$ of them occur), and the flag term $\sum_\sigma\langle Q_\sigma,M_\sigma(P)\rangle$ averages to zero within sampling error. So the construction puts almost all of its weight on patterns where the certificate has no slack.

\begin{proof}[Proof of Theorem~\ref{thm:main}]
Fix an atomless law $\mu$ and let $x_P$ be its four-point pattern
probabilities. Lemma~\ref{lem:psd}, applied to the colored flags, gives
$\sum_Px_P\langle Q_\sigma,M_\sigma(P)\rangle\ge0$ for each $\sigma$.
Averaging \eqref{eq:direct-four-cert} therefore yields
\[
 \Prob_{\mu^{\otimes4}}[X_1+X_2+X_3<2X_4]
 =\sum_Px_P\obj(P)
 \le\sum_Px_P\left(\obj(P)+
       \sum_\sigma\langle Q_\sigma,M_\sigma(P)\rangle\right)
 \le B.
\]
Atomless approximation as in Section~\ref{sec:example} extends the
bound to all laws on $\R_+$.
The construction of \citet{BF} gives $C\ge B$, hence $C=C_{\mathrm{BF}}$.
\end{proof}

The finite assertions \eqref{eq:direct-four-cert} are checked by
\path{anc/proofs/s3/0.400695372__c2/verify.py}: it replays the cover, reconstructs
the flag coefficients, and checks positive definiteness and every
pattern inequality in exact arithmetic over $\mathbb Q(r)$.
An independent verifier in \path{anc/proofs/s3/0.400695372__c2/independent/}
confirms them.
The entire proof, including the probability argument and the lower
construction, is formalized by \texttt{BeatAvg.Four.C\_eq\_Bstar} in
Lean~4. Its finite checks use the kernel, and its only axioms are
\texttt{propext}, \texttt{Classical.choice}, and \texttt{Quot.sound};
reproduction instructions are in \path{anc/lean/BeatAvg/Four/README.md}.
This proof uses only flag positivity and atomless approximation.
Section~\ref{sec:near-plateaus} uses near-optimality constraints for
the five additional exact $S_2$ plateaus in Table~\ref{tab:steps-bounds}.

%======================================================================
%======================================================================
\section{Flag SDP constraints from near-optimal laws}\label{sec:thmD}
\label{toplevel_sec:thmD}
%======================================================================
We call a flag certificate \emph{ordinary} if it uses only flag positivity,
so that it bounds every law, and \emph{augmented} if it also uses the
constraints of this section.
The ordinary colored flag certificates of Section~\ref{sec:pair-colors}
prove the first $S_2$ plateau and the exact four-roll value.
For the five additional exact $S_2$ plateaus in Table~\ref{tab:steps-bounds},
we use further information available when a law approaches the supremum.
Conditioning, moving mass, and replacing a lower tail give additional
constraints. We explain these properties, how to record their expectations
as linear SDP constraints, and how a fixed-point argument handles their
dependence on the unknown optimum.

\subsection{Properties of near-optimal laws}\label{sec:near-properties}
Fix $n\in\{3,4\}$ and $1<c\le2$, and write
\[
 P(\mu)=\Prob_{\mu^{\otimes n}}[X_1+\cdots+X_{n-1}<cX_n],
 \qquad S=\sup_\mu P(\mu).
\]
The laws are on $\R_+$ and may have atoms. The cases $(n,c)=(4,2)$
and $(3,c)$ give $S=C$ and $S=S_2(c)$, respectively.
For a sorted sample of $n$ distinct points, let $\sigma$ be the number
of winning target positions. Only the two largest positions can win,
so $P(\mu)=\E\sigma/n$ for atomless $\mu$ and $S\le2/n$.
The atomless approximation of Section~\ref{sec:example} lets us choose
atomless laws with $P(\mu_m)\to S$.

\subsubsection{Conditioning on a lower tail}
For $b>0$ put $F(b)=\mu[0,b)$ and
\[
 G(b)=\Prob[X_1+\cdots+X_{n-1}<cX_n,\ \max_iX_i<b].
\]
\begin{lemma}[Sub-law]\label{lem:sublaw4}
$G(b)\le S F(b)^n$.
\end{lemma}
\begin{proof}
If $F(b)>0$, conditioning on all draws lying below $b$ gives the iid
law $\mu(\,\cdot\mid[0,b))$, whose success probability is at most $S$.
If $F(b)=0$, both sides vanish.
\end{proof}

\subsubsection{Moving mass to a better point}
Let $d_\mu(z)$ be the success probability with one position, chosen
uniformly among the $n$ positions, fixed at $z$, and the other draws
independent with law $\mu$. Then $d_\mu$ is measurable,
$0\le d_\mu\le1$, and $\int d_\mu\,d\mu=P(\mu)$.
For atomless $\mu$, it is also the conditional mean of $\sigma/n$
given one point $z$.

\begin{lemma}\label{lem:degree4}
If $P(\mu)\ge S-\delta$ with $0\le\delta<1/(n-1)^2$, then
\[
 d_\mu(z)\le
 \frac{P(\mu)+\sqrt\delta/n}{1-(n-1)\sqrt\delta}
 \quad(z\ge0).
\]
\end{lemma}
\begin{proof}
For $\delta>0$ set $t=\sqrt\delta$ and
$\nu=(1-t)\mu+t\delta_z$. Keeping the terms with at most one draw
from $\delta_z$ and using Bernoulli's inequality gives
\[
 S\ge P(\nu)\ge(1-t)^nP(\mu)+nt(1-t)^{n-1}d_\mu(z)
 \ge(1-nt)P(\mu)+nt(1-(n-1)t)d_\mu(z).
\]
Rearrange using $S-P(\mu)\le t^2$. For $\delta=0$, apply the bound
with every sufficiently small positive $\delta'$ and let $\delta'\to0$.
\end{proof}

\begin{corollary}\label{cor:degree4}
Let $P(\mu_m)\to S$, let $L\le S$, and let
$p_m:\R_+\to[0,1]$ be measurable. Then
\[
 \liminf_m\int(d_{\mu_m}-L)p_m\,d\mu_m\ge0,\qquad
 \liminf_m\int(S-d_{\mu_m})p_m\,d\mu_m\ge0.
\]
The first assertion also holds with $L=S$.
\end{corollary}
\begin{proof}
Put $e_m=\sup_z d_{\mu_m}(z)-P(\mu_m)\ge0$.
Lemma~\ref{lem:degree4} gives $e_m\to0$. Since
$d_{\mu_m}-P(\mu_m)-e_m\le0$ and $0\le p_m\le1$,
\[
 \int(d_{\mu_m}-L)p_m\,d\mu_m
 \ge-e_m+(P(\mu_m)-L)\int p_m\,d\mu_m.
\]
Its lower limit is nonnegative. The second assertion follows from
$S-d_{\mu_m}\ge S-P(\mu_m)-e_m$.
\end{proof}

\subsubsection{Replacing a lower tail}
Define the crossing probability
\[
 \mathrm{BC}(b)=\Prob[X_1+\cdots+X_{n-1}<cX_n,\
                         X_n<b\le\max_{i<n}X_i].
\]
\begin{lemma}[Substitution]\label{lem:subst4}
$G(b)+\mathrm{BC}(b)\ge F(b)^nP(\mu)-(S-P(\mu))$.
\end{lemma}
\begin{proof}
The cases $F(b)=0,1$ are immediate. Otherwise replace the mass below
$b$ by a copy of $\mu$ scaled by $\varepsilon>0$, keeping the
conditional law above $b$ and its mass unchanged.
The new success probability is at most $S$. When all draws come
from the replaced part, its contribution is $F(b)^nP(\mu)$.
When the target comes from the upper part, sending $\varepsilon$
to zero preserves every original strict success in that part,
because the summands are nonnegative. Fatou's lemma therefore gives
\[
 S\ge F(b)^nP(\mu)+
       \Prob[X_1+\cdots+X_{n-1}<cX_n,\ X_n\ge b].
\]
Subtract this last probability from $P(\mu)$ to obtain the claim.
\end{proof}

\begin{lemma}[Crossings]\label{lem:bc4}
On a crossing success exactly one summand is at least $b$.
It is the largest point, and the target is the second largest.
For atomless $\mu$ and a sorted sample $t_1<\cdots<t_n$,
\[
 \mathrm{BC}(b)=\frac1n\Prob[
 t_1+\cdots+t_{n-2}+t_n<ct_{n-1},\
 t_{n-1}<b\le t_n].
\]
\end{lemma}
\begin{proof}
If the target is $x<b$ and a summand is $s\ge b$, the sum of the
other summands is less than $cx-s\le2x-s<x$.
Thus each of them is below $x$, whereas $s>x$.
The formula follows by exchangeability.
\end{proof}
Substitution and crossings, together with sub-law and degree constraints,
are used in the additional exact $S_2$ certificates below.

\subsection{Translating the properties into SDP constraints}\label{sec:near-sdp}
Continue with $n$, $c$, $S$ and $P(\mu)$ from
Section~\ref{sec:near-properties}. Fix a sample size $s\ge n+1$.
Record additive comparisons and pairwise ratio cutoffs as in
Section~\ref{sec:pair-colors-general}, so the pattern determines the
threshold and degree tests below. A certified finite cover $\cR$
contains every realizable pattern.

\subsubsection{Threshold tests}
Choose one root $z$ and $n$ further points, sorted as
$t_1<\cdots<t_n$. For each recorded cutoff, its inverse, and $q=1$,
define pattern functions by averaging over all such choices:
\[
 \begin{split}
 n_q(P)&=\operatorname{av}\mathbf1\{t_n<qz\},\\
 g_q(P)&=\operatorname{av}\frac{\sigma(t_1,\ldots,t_n)}n
                                      \mathbf1\{t_n<qz\},\\
 b_q(P)&=\operatorname{av}\frac1n
       \mathbf1\{t_1+\cdots+t_{n-2}+t_n<ct_{n-1},\
                         t_{n-1}<qz\le t_n\}.
 \end{split}
\]
The pair colors determine all comparisons with $qz$.
Conditioning on the root $Z\sim\mu$ gives
\[
 \E n_q=\E F(qZ)^n,\qquad
 \E g_q=\E G(qZ),\qquad
 \E b_q=\E\mathrm{BC}(qZ).
\]
The sub-law inequality therefore gives
$\E[S n_q-g_q]\ge0$ for every law.
For a known lower bound $L\le S$, the substitution inequality gives
\[
 \E[g_q+b_q-Ln_q]\ge-2(S-P(\mu)).
\]
This follows uniformly in the threshold from
$G+\mathrm{BC}-LF^n\ge F^n(P(\mu)-L)-(S-P(\mu))$.
Thus the substitution term has nonnegative lower limit along
near-optimal laws.

\subsubsection{Weighted degree tests}
Choose an ordered pair $(z,y)$ and further draws
$W_1,\ldots,W_{n-1}$. Its class $j$ records its orientation and ratio
color. Average over all these choices to define
\[
 k_j(P)=\operatorname{av}\mathbf1\{(z,y)\text{ has class }j\},\qquad
 d_j(P)=\operatorname{av}\frac{\sigma(z,W_1,\ldots,W_{n-1})}n
                         \mathbf1\{(z,y)\text{ has class }j\},
\]
where the arguments of $\sigma$ are sorted.
Writing $p_j(z)=\Prob[(z,Y)\text{ has class }j]$, independence gives
\[
 \E k_j=\int p_j\,d\mu,\qquad
 \E d_j=\int p_jd_\mu\,d\mu.
\]
Corollary~\ref{cor:degree4} shows that both $d_j-Lk_j$ and
$S k_j-d_j$ have nonnegative lower limits of their expectations
along near-optimal laws.

\subsubsection{The augmented certificate}
Choose nonnegative multipliers $a_q,b'_q,u_j,v_j$ and positive
semidefinite flag matrices $Q_\rho$. For a variable $x$ put
\begin{equation}\label{eq:plateau-near-score}
 \begin{split}
 \Lambda_P(x)={}&\obj(P)+\sum_\rho\langle Q_\rho,M_\rho(P)\rangle\\
 &+\sum_q\bigl[a_q(xn_q-g_q)+b'_q(g_q+b_q-Ln_q)\bigr](P)\\
 &+\sum_j\bigl[u_j(d_j-Lk_j)+v_j(xk_j-d_j)\bigr](P).
 \end{split}
\end{equation}
All functions of $P$ in this expression are fixed rational counts,
recomputed from the pattern definitions. Only the multipliers and
matrices are chosen by the SDP.
At $x=S$, every added term has nonnegative expectation or
nonnegative lower limit of its expectation.
Equivalently, if $x_P$ denotes pattern probabilities, these tests
give linear inequalities such as
$\sum_Px_P(Sn_q-g_q)(P)\ge0$, with vanishing errors for the
near-optimality tests. Their nonnegative linear combinations give
the augmented certificate.

Fix $L\le S\le\bar S$ and $x_0\le\bar S$.
For a fixed margin $\eta>0$, imposing
\[
 Q_\rho\succeq0,\qquad
 \Lambda_P(x_0)\le x_0,\qquad
 \Lambda_P(\bar S)\le\bar S-\eta\quad(P\in\cR),
\]
is an SDP feasibility problem. One may instead maximize $\eta\ge0$
and then verify that the resulting margin is positive.
The next subsection explains why these two finite sets of
inequalities prove $S\le x_0$.

\begin{remark}[Pattern restrictions from small slack]\label{sec:structure}
For an ordinary flag certificate with bound $\lambda$, write
$\mathrm{slack}(P)=\lambda-\val(P)\ge0$. Flag positivity gives
\begin{equation}\label{eq:slack}
 P(\mu)\le\lambda-\sum_Px_P(\mu)\,\mathrm{slack}(P).
\end{equation}
Consequently, the probability of a pattern with slack at least
$a>0$ is at most $(\lambda-P(\mu))/a$.
When $\lambda$ is close to the optimum, near-optimal laws therefore
concentrate on the small-slack patterns. This can guide the choice
of recorded relations and matrix kernels; the certificate is still
checked on the complete cover.
\end{remark}

\subsection{From near-optimality to an exact bound}\label{sec:idea}
Fix a sample size and a certified finite pattern cover $\cR$.
Add nonnegative multiples of the preceding constraints to the usual
flag squares, obtaining a pattern score $\Lambda_P(x)$ affine in
$x$, which stands for the unknown optimum $S$.
The multipliers and any known lower bound $L\le S$ are fixed before
$x$ varies. Define $\psi(x)=\max_{P\in\cR}\Lambda_P(x)$.

\begin{theorem}[Theorem D]\label{thm:D}
Suppose that along atomless laws $\mu_m$ with $P(\mu_m)\to S$,
every added term in $\Lambda_P(S)-\obj(P)$ has nonnegative
expectation or nonnegative lower limit of its expectation.
Then $S\le\psi(S)$.
\end{theorem}
\begin{proof}
The objective has expectation tending to $S$. There are finitely
many added terms, so
$S\le\liminf_m\E_{\mu_m}\Lambda_P(S)\le\psi(S)$.
\end{proof}

\begin{corollary}[Two-point criterion]\label{cor:twopoint4}
Under the hypotheses of Theorem~\ref{thm:D}, let $S\le\bar S$ and
$x\le\bar S$. If $\psi(x)\le x$ and $\psi(\bar S)<\bar S$,
then $S\le x$.
\end{corollary}
\begin{proof}
The function $\psi(u)-u$ is convex, nonnegative at $S$, and strictly
negative on $(x,\bar S]$ by the endpoint conditions.
Thus $S$ cannot belong to that interval.
\end{proof}
For exact $S_2$ plateaus we use $\bar S=2/3$ and $x=L$, the
recursive lower value. A strictly positive margin at $2/3$ excludes
every possible optimum above $L$.
For numerical bounds, if $\psi(L)\le L+e$ and
$\psi(\bar S)\le\bar S-\delta$ with $e\ge0$, $\delta>0$, the same
convexity argument gives
\[
 S\le L+\frac{e(\bar S-L)}{e+\delta}.
\]
One may also replace the known lower bound $L$ by the parameter $x$
in the substitution and lower degree terms. At $x=S$ their expectations
still have nonnegative lower limits, so the same criterion applies.
This is the full-strength version used for the narrow bands at
$\tau_0$ and $\tau_1$.

\subsection{Application to five exact plateaus}\label{sec:near-plateaus}
For $S=S_2(c)$, the constraints of Section~\ref{sec:near-sdp}
use four sample points: a root and three event draws, or an ordered
pair and two further draws.
For the five intervals below, Proposition~\ref{prop:onescale} gives
the uniform recursive lower value
\[
 L=V(\beta)=\frac23\left(\sqrt{7-15\beta+9\beta^2}-2+3\beta\right).
\]
The corresponding four-point certificates verify
$\psi(L)=L$ and $\psi(2/3)<2/3$ exactly.
\[
\begin{array}{c|c|r}
\text{interval}&\beta&\text{retained patterns}\\ \hline
(c_1,c_2]&4/27&179\\
(c_2,c_3]&\sqrt3/9&718\\
(\varphi,\alpha_1]&2\sqrt3/9&365\\
(\alpha_1,\alpha_2]&(216+38\sqrt{19})/961&1774\\
(\alpha_2,\alpha_3]&2s-3s^2+(1+\gamma)s^3&1995
\end{array}
\]
In the last row $\gamma=\sqrt3/9$ and
$s=(3-\sqrt{3-6\gamma})/(3(1+\gamma))$.
Corollary~\ref{cor:twopoint4} proves $S_2(b)\le L$ at each right
endpoint $b$. Monotonicity and the uniform lower construction give
$S_2(c)=L$ throughout the corresponding interval.
This proves the five additional zero-gap rows of Table~\ref{tab:steps-bounds}.

The standalone scripts
\path{verify_second_plateau.py},
\path{verify_third_plateau.py}, \path{verify_alpha1_plateau.py},
\path{verify_alpha2_plateau.py}, and \path{verify_alpha3_plateau.py}
check the pattern covers, flag coefficients,
nonnegative multipliers, positive semidefinite matrices, lower-value
identities, and both endpoint inequalities in exact algebraic arithmetic.
They are in \path{anc/proofs/s2/}.
The narrow bands at $\tau_0$ and $\tau_1$ use the same fixed-point
criterion with a nonzero first-endpoint error; they do not assert
exact plateau equality.

\section{Discontinuity points of $S_2(c)$ and $S_3(c)$}
\subsection{Two-against-one: The nested staircase structure}\label{sec:two}
%======================================================================

Throughout, $1<c<2$ unless stated otherwise. We derive the thresholds in
Table~\ref{tab:steps-thresholds} from the feasibility of explicit
lower-bound constructions, then combine those constructions with endpoint
certificates to obtain Table~\ref{tab:steps-bounds}. The two tasks are
distinct: A new winning pattern gives a lower bound immediately to the
right of a threshold; an upper bound at the threshold is needed to prove
a discontinuity.

For each law $\mu$, the probability $\Prob_\mu[X_1+X_2<cX_3]$ is
nondecreasing and left-continuous in $c$, since the strict event is the
increasing union of the events with $c'<c$. Taking suprema preserves
these properties. Thus an upper bound at $b$ holds for every $c\le b$,
and a lower construction valid for every $c>a$ gives a bound throughout
$(a,b]$. This is how we assemble the interval rows of
Table~\ref{tab:steps-bounds}.

\subsubsection{From one scale to a distribution}
\begin{figure}[ht]
\centering
\includegraphics[width=0.92\textwidth]{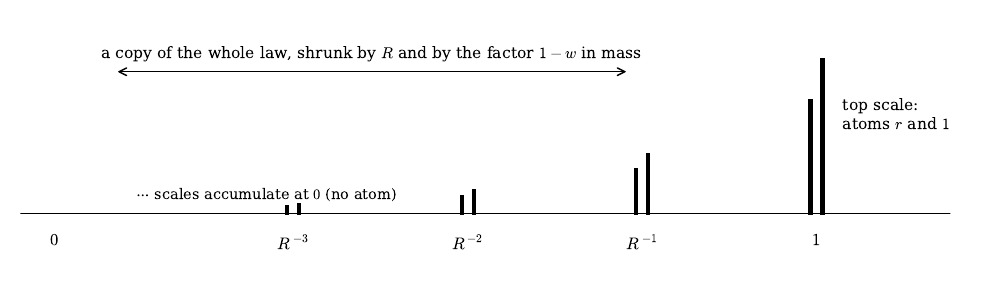}
\caption{A ladder of scales, the shape of our constructions for two against one when $c>\varphi$ (schematic, logarithmic horizontal axis): on each scale, two atoms at ratio $r\in(1/c,c-1)$, and each scale a copy of the scale above, shrunk by a large factor in value and by the factor $1-w$ in mass. The scales accumulate at $0$, which carries no mass.}
\label{fig:ladder}
\end{figure}

Our lower bounds for $1<c<2$ repeat a finite pattern of atoms on widely
separated scales (Figure~\ref{fig:ladder}). There are three choices:
The atom locations determine which triples win, their masses determine
the within-scale success probability, and the mass assigned to each
scale determines the effect of the recursion. Whether all optimal
values arise from such laws remains open.

 Proposition~3.3 of \citet{BF} reduces such laws to a single scale: for every law $\nu$ with finitely many atoms there is a law $\mu$ with
\begin{equation}\label{eq:boost}
  \Prob_\mu[X_1+X_2<cX_3]\;\ge\;\frac{\Prob_\nu[X_1+X_2<cX_3]}{1-\Prob_\nu[X_1+X_2=cX_3]}
\end{equation}
\cite[Proposition~3.3]{BF}; the law $\mu$ is the law of $\sum_j\eta^{j-1}Y_j$ with $Y_j$ iid with law $\nu$ and $\eta$ small, a construction going back to jlewk and W.~Sawin on \cite{MO474916}.

\begin{proposition}[One-scale reduction]\label{prop:onescale}
Let $\rho$ be a law with finitely many atoms in $[a,1]$, where $ca>1$, such that no atoms satisfy $y_1+y_2=cy_3$, and let $\beta=\Prob_\rho[Y_1+Y_2<cY_3]$. Then
\[
  S_2(c)\;\ge\;V(\beta)\;:=\;\sup_{w\in(0,1)}\frac{1-(1-\beta)w^2}{3-3w+w^2}\;=\;\tfrac23\Bigl(\sqrt{9k^2-3k+1}-3k+1\Bigr),\qquad k=1-\beta .
\]
\end{proposition}

\begin{proof}
Let $\nu=(1-w)\delta_0+w\rho$. If $X_3=0$ the event fails. If $X_3$ is an atom of $\rho$, a pair containing a zero wins, since the other point is at most $1<cX_3$, and a pair of atoms of $\rho$ wins with conditional probability $\beta$. Ties occur only when all three points are $0$. So $\Prob_\nu[\text{win}]=w(1-w^2+w^2\beta)$ and $\Prob_\nu[\text{tie}]=(1-w)^3$, and \eqref{eq:boost} gives $S_2(c)\ge w(1-(1-\beta)w^2)/(1-(1-w)^3)=(1-(1-\beta)w^2)/(3-3w+w^2)$ for every $w$. The closed form is a calculus exercise.
\end{proof}

$V$ is increasing in $\beta$. We can therefore count any guaranteed
collection of winning triples, obtain a lower bound on $\beta$, and
apply $V$ to it. The locations may depend on $c$ while the masses
remain fixed. Once the same comparisons are feasible for every $c>a$,
the resulting lower bound is uniform to the right of $a$.

\begin{corollary}\label{cor:steps}\leavevmode
\begin{enumerate}
  \item For every $c\in(1,2)$: $S_2(c)\ge V(0)=(2\sqrt7-4)/3=0.4305008\ldots$ (take $\rho=\delta_1$).
  \item For $\sqrt2<c\le\varphi$: $S_2(c)\ge V(4/27)=(2\sqrt{403}-28)/27=0.4499896\ldots$.
  \item For $c>c_2$: $S_2(c)\ge V\bigl(1/(3\sqrt3)\bigr)=0.4573615\ldots$.
  \item For $c>\varphi$: $S_2(c)\ge V\bigl(2\sqrt3/9\bigr)=0.5030984\ldots$.
\end{enumerate}
\end{corollary}

\begin{proof}
(1) With one atom no pair of atoms wins, so $\beta=0$. This is the recipe of \cite[Section~4]{BF} applied to a Bernoulli law, as Yu did for three against one \cite{Yu}. The value is also attained by an actual law: the geometric law $\Prob[X=R^{-d}]=w(1-w)^d$, $d\ge0$, with $R>\max\{2,1/(c-1)\}$ and the optimal $w=(4-\sqrt7)/3$.

(2) Take $\rho=\frac23\delta_r+\frac13\delta_1$ with $1/c<r<c/2$, possible exactly when $c>\sqrt2$, and $r$ avoiding the finitely many tie values. Against $Y_3=1$ the pair $(r,r)$ wins since $2r<c$, the pair $(r,1)$ loses since $1+r>1+1/c\ge c$ for $c\le\varphi$, and $(1,1)$ loses; against $Y_3=r$ every pair loses, since $2r>cr$ and $1+r>cr$. So $\beta=(2/3)^2\cdot\frac13=4/27$.

(3) Take $\rho=p\,\delta_a+q\,\delta_b+q\,\delta_1$ with $p=1/\sqrt3$, $q=(1-p)/2$ and $1/c<a$, $2a<cb$, $a+b<c$, possible exactly when $1/c+2/c^2<c$, that is $c^3>c+2$. The pairs $(a,a)$ against $b$ and against $1$, and $(a,b)$ against $1$, win, so $\beta\ge p^2q+p^2q+2pq^2=2pq(p+q)=1/(3\sqrt3)$.

(4) Take $\rho=(1-p)\delta_r+p\,\delta_1$ with $p=1/\sqrt3$ and $1/c<r<c-1$, possible exactly when $c^2>c+1$. Against $Y_3=1$ the pairs $(r,r)$ and $(r,1)$ win and $(1,1)$ loses; against $Y_3=r$ every pair loses. So $\beta=p\bigl((1-p)^2+2p(1-p)\bigr)=2\sqrt3/9$. For $c\ge2$ the bound follows by monotonicity.
\end{proof}

\subsubsection{A complete example: The interval $(c_1,c_2]$}
The second construction in Corollary~\ref{cor:steps} illustrates the
whole proof. Its positive block has atoms $r,1$ with masses $2/3,1/3$.
To make the pair $(r,r)$ win against $1$, while a zero paired with $1$
wins against $r$, we need
\[
  2r<c,\qquad 1<cr,
  \quad\text{or equivalently}\quad \frac1c<r<\frac c2.
\]
This window is nonempty exactly when $c>\sqrt2=c_1$.
For $c\le\varphi$ these are the only winning triples within the positive
block, so $\beta=4/27$ and the recursive lower value is
\[
  L=V(4/27)=\frac{2\sqrt{403}-28}{27}.
\]
It is a lower bound for every $c>c_1$: The construction gives it up to
$\varphi$, and monotonicity extends it beyond.

To prove equality on $(c_1,c_2]$, we need an upper bound at $c_2$.
The four-point colored certificate of Section~\ref{sec:near-plateaus}
uses the near-optimal-law constraints to verify
\[
  \psi(L)=L,\qquad \psi(2/3)<2/3.
\]
The two-point criterion gives $S_2(c_2)\le L$. Monotonicity then gives
$S_2(c)\le L$ for every $c\le c_2$, proving equality throughout
$(c_1,c_2]$. The script \path{verify_second_plateau.py} checks this
certificate in exact arithmetic. The construction determines the left
threshold $c_1$; the certificate at $c_2$ proves that its value remains
optimal through the right endpoint.

The next construction in Corollary~\ref{cor:steps} is feasible for every
$c>c_2$ and gives $V(\sqrt3/9)>L$. Hence it also proves a right-hand
jump at $c_2$. The same sequence of arguments underlies the other rows:
Feasibility, a uniform recursive lower bound, and a certified endpoint
upper bound.

\subsubsection{Three windows: $\sqrt2$, $c_2$ and $\varphi$}\label{sec:windows}
\begin{figure}[t]
\centering
\includegraphics[width=\textwidth]{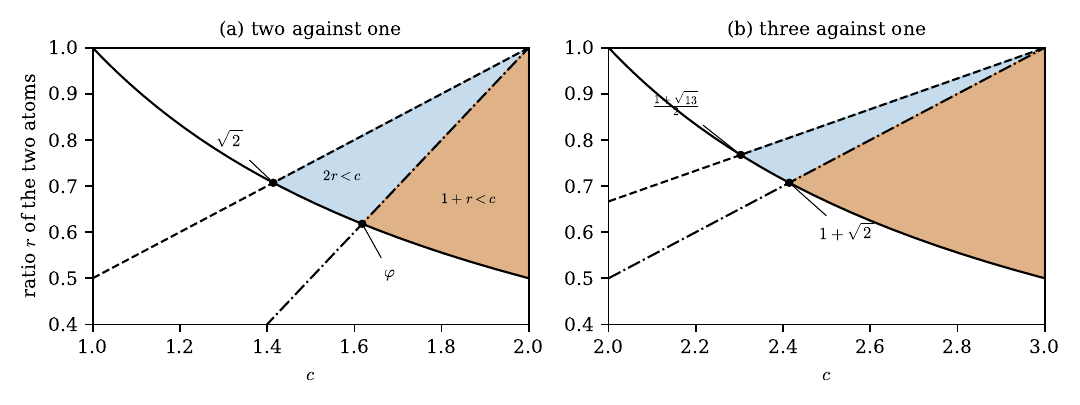}
\caption{Where the steps come from. (a) Two against one: the lines $r=1/c$ (solid), $r=c/2$ (dashed) and $r=c-1$ (dash-dot); the window $1/c<r<c/2$ (blue and tan) opens at $\sqrt2$, and its part $1/c<r<c-1$ (tan) at $\varphi$. (b) Three against one: $r=1/(c-1)$, $r=c/3$ and $r=(c-1)/2$; the windows open at $(1+\sqrt{13})/2$ and $1+\sqrt2$. Inside each window the constructions of Corollary~\ref{cor:steps} and Proposition~\ref{prop:three-above} give a bound that does not depend on $c$.}
\label{fig:windows}
\end{figure}

In the one-scale reduction, an atom $r<1$ must beat the top atom together with something negligible, $1<cr$, that is $r>1/c$; and it must help the top atom win, either because two copies of $r$ lose to it ($2r<c$) or because $r$ and the top atom together lose to it ($1+r<c$). The first window $1/c<r<c/2$ opens when $c^2>2$, the second, $1/c<r<c-1$, when $c^2>c+1$. With three atoms $a<b<1$ the condition $2a<cb$ adds a middle step, and the window opens when $c^3>c+2$. The window conditions guarantee the winning configurations used in our lower bounds, so within a window these constructions give a bound that does not depend on $c$; further configurations may become available as $c$ increases (for example, above $\varphi$ an atom in the first window may also satisfy $1+r<c$). For the two-atom constructions, the endpoints of the feasible interval meet at the threshold (Figure~\ref{fig:windows}). At $c_2$ a third atom provides an intermediate target; at $\varphi$ no extra atom is needed, since the existing pair $(r,1)$ becomes a winning pair against $1$. The transition concerns the winning pattern, not necessarily the number of atoms.

For $k$ against one, $\sup_\mu\Prob[X_1+\dots+X_k<cX_{k+1}]$ with $k-1<c<k$, the analogous second window, $(k-1)r+1<c$ together with $(c-k+2)r>1$, opens at the positive root of $c^2-(k-1)c-1=0$: the golden ratio for $k=2$, the silver ratio $1+\sqrt2$ for $k=3$, and $(3+\sqrt{13})/2$ for $k=4$.

\subsubsection{Longer ladders below the golden ratio}
The examples at $c_1$ and $c_2$ extend to a family of ladders. Before normalizing the largest atom to $1$, take positive atoms
\[
 a_0=1<a_1<\cdots<a_n,
 \qquad a_n<c,\qquad 1+a_j<ca_{j+1}\quad(0\le j<n).
\]
The first inequality ensures that every positive atom beats any positive atom paired with zero. The other inequalities make the smallest positive atom, paired with an atom at one level, lose to the next level. The limiting positions satisfy
\[
 r_0=1,\qquad r_{j+1}=\frac{1+r_j}{c},\qquad
 r_n=\frac{2+c+\cdots+c^{n-1}}{c^n}\quad(n\ge1).
\]
Every strict ladder has $a_n>r_n$. Conversely, when $r_n<c$, small positive perturbations of the recurrence give a strict ladder with $a_n<c$. Thus the feasibility threshold is precisely
\[
 c_n^{n+1}=2+c_n+\cdots+c_n^{n-1}.
\]
For $n=1$ this reads $c_1^2=2$, with the sum after $2$ empty. The next three equations are those for $c_2,c_3,c_4$ in Theorem~\ref{thm:steps}. The roots increase to $\varphi$: the recurrence has fixed point $1/(c-1)$, and this equals the allowed upper endpoint $c$ exactly when $c^2=c+1$. Equivalently, the threshold equation gives
\[
 c_n^n(c_n^2-c_n-1)=c_n-2.
\]

Choose masses on the ladder and zero, and recursively replace zero as in \eqref{eq:boost}. The guaranteed winning triples give a lower bound valid for every $c>c_n$, with the same masses but atom locations and scale separation depending on $c$. Additional comparisons can improve this bound: at $c_4$, the certified construction also uses $2a_1<ca_4$, which is compatible with the ladder just above the threshold. The resulting bounds above $c_3$ and $c_4$ are $0.4606036$ and $0.4627479$, respectively. These are finite constructions; no limiting ladder is needed for these two jumps.

More explicitly, put masses $p_0,\ldots,p_n$ summing to one on the
positive block. The ladder comparisons guarantee within-scale probability
\[
 B=p_0^2\sum_{k=1}^n p_k
   +2p_0\sum_{1\le i<k\le n}p_i p_k.
\]
For $n=4,5,6$, the comparisons $2a_1<ca_k$ for $k\ge4$ add
$p_1^2\sum_{k=4}^n p_k$. The certificates supply rational masses and a
scale weight $w$; evaluating
$[1-(1-B)w^2]/[3-3w+w^2]$ exactly gives the lower bounds immediately
after $c_n$. Some rows use the larger algebraic value obtained by
optimizing these masses and $w$, as in the worked example.

\subsubsection{Branching ladders below the golden ratio}
The extra comparison involving two copies of $a_1$ can also be inserted
earlier in the ladder. Alongside $r_0=1$, $r_{k+1}=(1+r_k)/c$, define
\[
 b_0=4/c^2,\qquad b_{j+1}=(1+b_j)/c.
\]
The positive block $1,r_1,r_2,b_0,\ldots,b_j$ introduces a target for
the pair $(r_1,r_1)$. Its last atom must lie below $c$, so its limiting
feasibility condition is $b_j<c$. Clearing denominators gives
\[
 \tau_j^{j+3}=4+\sum_{k=2}^{j+1}\tau_j^k\qquad(j\ge0).
\]
The empty sum gives $\tau_0=\sqrt[3]4$. For $1<c<\varphi$,
\[
 b_0-r_3=\frac{(2-c)(c-1)}{c^3}>0,\qquad
 r_4-b_0=\frac{(c-2)(c^2-c-1)}{c^4}>0.
\]
Iterating the same increasing map gives $r_{j+3}<b_j<r_{j+4}$.
Consequently
\[
 c_{j+3}<\tau_j<c_{j+4},\qquad \tau_j\longrightarrow\varphi.
\]
This explains the interlacing of the $c_n$ and $\tau_j$ rows in
Table~\ref{tab:steps-thresholds}.

For $\tau_j<c<\varphi$, make the comparisons strict by setting
\[
 x_0=1,\quad x_1=2/c+\varepsilon,\quad
 x_2=(1+x_1)/c+\varepsilon,\quad
 x_3=2x_1/c+\varepsilon,\quad
 x_{k+1}=(1+x_k)/c+\varepsilon\quad(k\ge3).
\]
For sufficiently small $\varepsilon>0$, the points are strictly ordered
and $x_{j+3}<c$. They satisfy $1+x_i<cx_k$ whenever $i<k$ and
$2x_1<cx_k$ for $k\ge3$. In particular, the comparison
$1+x_2<cx_3$ has limiting slack $(2-c)(c-1)/c^2>0$.
The limiting last position is $b_j$, so these requirements are feasible
exactly for $c>\tau_j$ in this range. With positive masses
$p_0,\ldots,p_{j+3}$ summing to one, the guaranteed probability is
\[
 B=p_0^2\sum_{k=1}^{j+3}p_k
   +2p_0\sum_{1\le i<k\le j+3}p_i p_k
   +p_1^2\sum_{k=3}^{j+3}p_k.
\]
Normalize the largest atom to $1$ and apply
Proposition~\ref{prop:onescale}. Fixed masses and a fixed scale weight
give the uniform lower bounds after $\tau_0,\ldots,\tau_3$;
monotonicity extends them to every $c>\tau_j$.
The extra summand in $B$ explains the improvement over the unbranched
ladder. The corresponding endpoint certificates prove four jumps;
feasibility alone does not prove a jump at every $\tau_j$.

\subsubsection{A second ladder above the golden ratio}
Above $\varphi$, the window $1/c<r<c-1$ permits the pair $(r,1)$ to win against $1$. We can insert a ladder into this window. Use zero, the largest atom $1$, and atoms
\[
 0<a_0<\cdots<a_n<1,\qquad
 ca_0>1,\qquad ca_{j+1}>a_0+a_j\ (0\le j<n),\qquad a_n<c-1.
\]
The last inequality ensures that every pair of positive atoms except $(1,1)$ wins against $1$; the middle inequalities add successes against the intermediate targets. The limiting positions now obey
\[
 r_0=1/c,\qquad r_{j+1}=\frac{r_0+r_j}{c},\qquad
 r_n=\frac{2+c+\cdots+c^{n-1}}{c^{n+1}}\quad(n\ge1).
\]
As before, the strict inequalities are feasible exactly when the limiting last position lies below the upper endpoint, here $r_n<c-1$. The resulting thresholds satisfy
\[
 \alpha_n^{n+2}-\alpha_n^{n+1}
 =2+\alpha_n+\cdots+\alpha_n^{n-1}\qquad(n\ge1).
\]
The first two are $\alpha_1^3=\alpha_1^2+2$ and
$\alpha_2^4=\alpha_2^3+\alpha_2+2$. These are the same thresholds obtained
from an unnormalized block
\[
 1=A_0<A_1<\cdots<A_n<H<c,\qquad
 1+A_i<cA_{i+1},\qquad A_n<(c-1)H.
\]
The limiting positions are $A_i=cr_i$ and $H=A_n/(c-1)$.
This form makes the mass calculation explicit. Give the lower ladder
total mass $s$ and the top atom $H$ mass $1-s$. If the conditional
success probability guaranteed within the lower ladder is $\gamma$,
then the entire positive block has guaranteed probability
\[
 \beta=\gamma s^3+(1-s)-(1-s)^3
       =2s-3s^2+(1+\gamma)s^3.
\]
Indeed, every pair with at least one member below $H$ wins against $H$.
Optimizing over $s$ gives
\[
 s=\frac{3-\sqrt{3-6\gamma}}{3(1+\gamma)}.
\]
Within the lower ladder, put conditional mass $t$ at $A_0$ and
$(1-t)/n$ at each of $A_1,\ldots,A_n$. Its guaranteed probability is
\[
 \gamma_n(t)=t^2(1-t)+\frac{n-1}{n}t(1-t)^2.
\]
For $n=1,2,3,4$, maximizing this polynomial gives, respectively,
\[
 \gamma=\frac4{27},\qquad \frac{\sqrt3}{9},\qquad
 \frac{-20+14\sqrt7}{81},\qquad
 \frac{-35+13\sqrt{13}}{54}.
\]
The maximizing $t$ values are $2/3$, $1/\sqrt3$,
$(\sqrt7-1)/3$ and $(\sqrt{13}-2)/3$.
Applying $V$ to these $\beta$ values gives the uniform lower bounds
after $\alpha_1,\ldots,\alpha_4$ in Table~\ref{tab:steps-bounds}.
The first two match their right-endpoint certificates exactly.
The lower bound after $\alpha_3$ is
$0.50905413288916\ldots$; its certificate at $\alpha_4$ leaves a
positive gap. Longer ladders give further constructive thresholds,
but we do not claim a discontinuity at every one of them.

The roots $\alpha_n$ increase to $\beta_3'$, the root of $c(c-1)^2=1$. Indeed,
\[
 \alpha_n^n\bigl[\alpha_n(\alpha_n-1)^2-1\bigr]=\alpha_n-2,
\]
and the fixed point of the limiting recurrence is $1/[c(c-1)]$, which meets $c-1$ at $\beta_3'$. At this accumulation point a different finite pattern becomes available. With atoms $0,a,b,1$, require
\[
 1/c<a,\qquad a+b<cb,\qquad b<c-1.
\]
Now $b$ itself, paired with $a$, loses to the target $b$; the ladder can close with this comparison instead of requiring another, larger target. The conditions are feasible exactly when
\[
 \frac1c<(c-1)^2,
\]
since they require $1/c<a<(c-1)b<(c-1)^2$, and these inequalities also allow a choice of $a,b$. This finite family gives $S_2(c)\ge0.5162682$ for every $c>\beta_3'$.

\subsubsection{The larger thresholds and the nested structure}
The thresholds $1,\varphi,\beta_3$ belong to another family. Take atoms
\[
 x_0=0<x_1<\cdots<x_m=1
\]
and require
\[
 1+x_{k-1}<cx_k\qquad(1\le k\le m).
\]
Then a pair wins against $x_k$ exactly when at least one member of the pair has index smaller than $k$. Indeed, such a pair sums to at most $1+x_{k-1}$; otherwise its sum is at least $2x_k>cx_k$. Iteration forces
\[
 x_{m-1}>\sum_{j=1}^{m-1}c^{-j},\qquad x_{m-1}<c-1.
\]
For $m\ge2$, these conditions are feasible exactly when
\[
 c^m>1+c+\cdots+c^{m-1}.
\]
For sufficiency, choose the intermediate atoms slightly above the partial sums of $\sum_{j\ge1}c^{-j}$, preserving the strict final inequality. Thus the roots $\beta_m$ of the corresponding equalities are feasibility thresholds, with $\beta_1=1$, $\beta_2=\varphi$ and $\beta_3$ as in Theorem~\ref{thm:steps}. The case $m=1$ is the law on $\{0,1\}$ used for the first plateau. If $p_k$ are the masses and $T_k=\sum_{j=k}^m p_j$, the base success probability is
\[
 F=\sum_{k=1}^m p_k(1-T_k^2),
\]
and recursive zero replacement gives the lower bound $F/(1-p_0^3)$. For $m=3$, rational masses give $S_2(c)\ge0.5441582$ for every $c>\beta_3$.

This gives all the threshold families in Table~\ref{tab:steps-thresholds}:
The $c_n$ and $\tau_j$ accumulate at $\beta_2=\varphi$, the $\alpha_n$
accumulate at $\beta_3'<\beta_3$, and the $\beta_m$ increase to $2$.
Theorem~\ref{thm:steps} proves jumps at the eighteen listed points.
Whether every member of these families is a discontinuity remains open.

\subsubsection{Certifying the right endpoints}\label{sec:colored-jumps}
An endpoint certificate must hold at the algebraic threshold itself.
An upper bound at a nearby coefficient below it cannot be extended
to the threshold by monotonicity. The verifiers therefore work in the
ordered algebraic field of the endpoint, representing its root by a
defining polynomial and a rational isolating interval.

\paragraph{Ordinary colored certificates.}
Record both the additive triple comparisons and pairwise ratio colors
as in Section~\ref{sec:pair-colors-general}. These relations subdivide
the possible samples into finitely many patterns. The cutoffs are
chosen from the ratios appearing in the limiting constructions:
They let the flags distinguish positions that the additive comparisons
alone do not distinguish. For the four-point certificates at the
branching endpoints, write $r_0=1$, $r_{k+1}=(1+r_k)/c$ and
$b_0=4/c^2$, $b_{k+1}=(1+b_k)/c$, and use the sorted distinct cutoffs
\[
 \left\{\frac{c+2}{2c},r_1,\ldots,r_{j+2},b_0,\ldots,b_{j-1},c\right\}
 \quad\text{at }c=\tau_j,
\]
where the list of $b$'s is empty for $j=0$, and
\[
 \left\{\frac{c+2}{2c},r_1,\ldots,r_n,b_0,\ldots,b_{n-4},c\right\}
 \quad\text{at }c=c_n\ (n=4,5,6).
\]
At $c=\alpha_3$ the cutoffs are
\[
 \left\{2/c,c,1/(c-1),(c+2)/(2c),(c+2)/c^2,
              (c^2+c+2)/(c(c+2))\right\}.
\]
The flag families have root and flag sizes $(t,f)=(0,2),(1,2),(2,3)$.
These choices need not be minimal. What matters for the proof is a
complete cover of the patterns and positive semidefinite matrices.
Every excluded sign prefix carries a Gordan certificate checked in
the endpoint field; every retained pattern is tested against the score
bound. The matrices are supplied as rational Gram factors, so their
positive semidefiniteness follows without a numerical eigenvalue test.
The verifier reconstructs the flag coefficients and checks
\[
 \operatorname{obj}(P)+\sum_\sigma
 \langle Q_\sigma,M_\sigma(P)\rangle\le U,
 \qquad Q_\sigma=L_\sigma L_\sigma^\top/10^{18},
\]
on every retained pattern, with integer matrices $L_\sigma$.
Section~\ref{sec:pair-colors-general} then gives $S_2(c)\le U$.
This argument uses constraints valid for every law. It gives the
four-point upper bounds at $c_4,c_5,c_6,\tau_2,\tau_3$ used in
Table~\ref{tab:steps-bounds}, and also valid bounds at $\tau_0,\tau_1$
and $\alpha_3$ that the augmented certificates below improve.

At $\varphi,\beta_3'$ and $\beta_3$, Table~\ref{tab:steps-bounds}
uses six-point certificates that record only additive comparisons.
Their matrices also have rational Gram factors, while pattern
completeness is checked in the exact endpoint field. The corresponding
data and verifiers are in \path{anc/proofs/s2/six_point_endpoints/}.
For a small illustration of algebraic completeness,
Appendix~\ref{app:gordan-phi} gives the four-point additive cover at
$\varphi$: Of its $32$ exclusion certificates, only two require
weights in $\mathbb Q(\sqrt5)$.

The four-point proof data and lower constructions are in
\path{anc/proofs/s2/}. The standard-library scripts
\path{verify_cuberoot_discontinuity.py},
\path{verify_family_discontinuities.py}, and
\path{verify_ladder_discontinuities.py} check the covers, upper
certificates, lower probabilities, and positive differences exactly.

\paragraph{Certificates using near-optimal laws.}
For the five further exact plateaus, the same four-point framework
also records the constraints of Section~\ref{sec:near-sdp}.
They hold along laws approaching the supremum, and their score
$\Lambda_P(x)$ depends affinely on the unknown optimum $x$.
With $\psi(x)=\max_P\Lambda_P(x)$, Theorem~\ref{thm:D} gives
$S_2(c)\le\psi(S_2(c))$. An ordinary bound on every pattern at one
fixed parameter would not by itself resolve this dependence.

Let $L$ be the uniform lower value from the construction. The
certificates in Section~\ref{sec:near-plateaus} verify
$\psi(L)=L$ and $\psi(2/3)<2/3$ exactly. Since
$\psi(x)-x$ is convex, it is strictly negative for
$L<x\le2/3$. The general bound $S_2(c)\le2/3$ therefore forces
$S_2(c)\le L$. This proves the endpoint equality needed for an
exact plateau. The scripts listed in that section check the algebraic
lower values, nonnegative multipliers, PSD matrices, pattern covers,
and both score inequalities.

The bands ending at $\tau_0,\tau_1$ and $\alpha_4$ use the same
argument with a positive first-endpoint error. If
\[
 \psi(L)\le L+e,\qquad
 \psi(2/3)\le2/3-\delta,\qquad e\ge0,\quad\delta>0,
\]
then
\[
 S_2(c)\le L+\frac{e(2/3-L)}{e+\delta}.
\]
These inequalities, rather than the residuals of a numerical SDP
solve, certify the reported gaps. The $\tau_0$ certificate uses the
full-strength substitution and degree constraints of
Section~\ref{sec:idea}; \path{verify_near_small.py} checks
\path{near_tau0_refined_cover.json} with
\path{near_tau0_full_near_certificate.json}.
For $\tau_1$, \path{verify_tau1_band.py} with
\texttt{--ratios 31} \texttt{--full-\allowbreak near-\allowbreak optimality}
checks the bound in the table.
The script \path{verify_alpha4_band.py} checks the band ending at
$\alpha_4$, and \path{verify_alpha4_discontinuity.py} checks its
separation from the construction available to its right.
The widths are less than $4.63\cdot10^{-9}$,
$1.24\cdot10^{-6}$ and $2.33\cdot10^{-7}$, respectively.
None of these three positive gaps establishes exact equality.

\subsubsection{Assembling the intervals and jumps}
Write $L_a$ for a lower value valid for every $c>a$, and $U_b$ for a
certified upper bound at $b$. For consecutive thresholds $a<b$,
monotonicity gives
\[
  L_a\le S_2(c)\le U_b\qquad(a<c\le b).
\]
Equality $L_a=U_b$ proves an exact plateau. A positive difference
$U_b-L_a$ is only a band; it does not prove that $S_2$ is constant.
At its right endpoint we compare with the next construction:
\[
 S_2(b+)-S_2(b)\ge L_b-U_b.
\]
This difference must be strictly positive to certify a jump. The
construction is uniform for all $c>b$, although its support and scale
separation may depend on $c$. A witness at one nearby rational
coefficient would not suffice.

The first row of Table~\ref{tab:steps-bounds} is elementary.
For $0<c\le1$, at most one of the three choices of target can win,
giving $S_2(c)\le1/3$. For the reverse bound, use the uniform law on
$N$ positive atoms separated by a ratio $R>\max\{1,2/c\}$.
Whenever the three draws have a unique largest atom, that atom beats
the sum of the other two after multiplication by $c$.
The probability of a repeated draw tends to zero as $N\to\infty$,
so the success probability tends to $1/3$.
The row $(1,c_1]$ combines $V(0)$ with the ordinary three-point
certificate of Section~\ref{sec:plateau}.
The five further zero-gap rows,
$(c_1,c_2]$, $(c_2,c_3]$, $(\varphi,\alpha_1]$,
$(\alpha_1,\alpha_2]$ and $(\alpha_2,\alpha_3]$,
combine the corresponding one-scale values with the near-optimal-law
certificates of Section~\ref{sec:near-plateaus}.

Below $\varphi$, the remaining rows alternate between the unbranched
$c_n$ construction and the $\tau_j$ construction. Ordinary four-point
certificates give the upper bounds at $c_4,c_5,c_6,\tau_2,\tau_3$;
the sharper bounds at $\tau_0$ and $\tau_1$ use the fixed-point
criterion. The bound at $\varphi$ is an ordinary six-point certificate.
Above $\varphi$, the $\alpha_3$ lower construction and the near-optimal
certificate at $\alpha_4$ give the narrow band
$(\alpha_3,\alpha_4]$. It is not an exact plateau.
The next row uses the $\alpha_4$ lower construction and the ordinary
six-point upper certificate at $\beta_3'$.
The final row uses the closing construction at $\beta_3'$ and the
ordinary six-point certificate at $\beta_3$.
The construction available for every $c>\beta_3$ supplies the last
jump bound.

All comparisons use unrounded values. Only afterwards do we round
lower bounds down, upper bounds up, and the displayed gaps and jump
bounds as specified in Table~\ref{tab:steps-bounds}.
Thus the eighteen threshold equations describe feasible patterns,
while the eighteen positive endpoint separations prove
Theorem~\ref{thm:steps}. Whether the full function follows the nested
staircase of these constructions remains open.

\begin{proof}[Proof of Theorem~\ref{cor:fair2}]
Table~\ref{tab:steps-bounds} gives $S_2(\varphi)\le0.4702903$;
monotonicity extends this to every $0<c\le\varphi$.
Section~\ref{sec:near-plateaus} proves the exact value
$V(2\sqrt3/9)=L_\varphi$ on $(\varphi,\alpha_1]$,
and monotonicity extends its lower bound to every $c>\varphi$.
The plateau gives the right limit, and
$L_\varphi-0.4702903>0.0328081$ gives the jump bound.
The two ranges place $S_2$ strictly below and above $1/2$,
respectively, proving the threshold and the absence of a fair coefficient.
\end{proof}

%======================================================================
\subsection{Three-against-one: Exact intervals and bounds}\label{sec:three}
%======================================================================

We now consider
\[
 S_3(c)=\sup_\mu\Prob[X_1+X_2+X_3<cX_4],\qquad c>0.
\]
As for $S_2$, the function is nondecreasing and left-continuous.
Figure~\ref{fig:three-against-one} displays four exact intervals in
green and the exact value at $c=2$ as a green point. Elsewhere,
the orange and blue curves are lower and upper bounds. A horizontal
segment of either bound need not be a plateau of $S_3$.
We explain the constructions, their feasibility thresholds, and the
endpoint certificates that produce these curves.
All upper certificates in this subsection use ordinary flag PSD
constraints, valid for every law.

An upper bound at $b$ extends to every $c\le b$, while a lower
construction valid for every $c>a$ extends to the right of $a$.
Matching bounds prove an exact interval $(a,b]$; bounds that do not
match leave the shaded uncertainty in the figure. The upper panel
shows the full range $0.9\le c\le2.6$, and the lower panel magnifies
$2\le c\le2.44$. The diagonal connectors identify this region.
The dashed horizontal line marks $1/2$: We conclude by explaining
why $S_3$ jumps across it at the silver ratio, answering
\cite[Problem~5.3]{BF}.

\subsubsection{Below two: Exact intervals and the remaining band}
For $0<c\le1$, at most one of four draws can beat the sum of the other
three after multiplication by $c$, so $S_3(c)\le1/4$.
Uniform laws on $N$ positive atoms separated by a ratio greater than
$3/c$ approach $1/4$ as $N\to\infty$, proving equality.
For $c>1$, a one-atom-per-scale construction gives
\begin{equation}\label{eq:s3-below-two}
 L_0=\max_{0<z<1}\frac{z^2(3-2z)}{1+z+z^2+z^3}
     =0.309684795417425\ldots,\qquad
 S_3(c)=L_0\quad(1<c\le\sqrt2).
\end{equation}
The maximizing $z=0.7139293479023\ldots$ is the unique root in $(0,1)$
of $5z^3+4z^2+3z-6=0$.

\begin{proof}
Take $\Prob[X=\rho^k]=(1-z)z^k$, $k\ge0$, with
$0<\rho<\min\{(c-1)/2,1/c\}$.
For $1<c\le2$, success occurs precisely when the target is at the
largest occupied scale and at most one summand is at that scale.
Thus
\[
 P=(1-z)z^3+3(1-z)^2z^2+z^4P.
\]
Solving for $P$ gives the displayed value. Its derivative has
numerator $-z(5z^3+4z^2+3z-6)$; the cubic increases strictly from
$-6$ to $6$ on $[0,1]$, so its unique zero maximizes this family.

For the matching upper bound at $\sqrt2$, record the pair comparisons
$x_j<\sqrt2 x_i$ on four sorted points, together with the additive
target comparisons. An exact certificate with two rooted PSD blocks
has every pattern score at most $L_0$. Ordinary conditional-square
positivity gives $S_3(\sqrt2)\le L_0$, and monotonicity proves equality
throughout $(1,\sqrt2]$. Appendix~\ref{app:s3-certificates} gives the
matrices and verification details.
\end{proof}

The endpoint $\sqrt2$ is sharp. For $\sqrt2<c<2$, the window
$1/c<r<c/2$ permits masses $(z,a,b)$ on $(0,r,1)$ with guaranteed
success numerator
\[
 (a+b)z^3+3(a+b)^2z^2+3za^2b.
\]
Recursive zero replacement divides this by $1-z^4$.
The rational masses in \path{below_two_uniform_jump.json} give
$S_3(c)>0.321369982677$ throughout this window; monotonicity extends
the bound to every $c>\sqrt2$.
It exceeds $L_0$ by more than $0.011685187$, proving the jump shown
after the second green interval.

Between $\sqrt2$ and $2$, the orange curve uses finite versions of
the ladder of \citet{BF}, with optimized rational masses and, for
some witnesses, an additional positive atom. The ladders satisfy
\[
 a_0=0<a_1<\cdots<a_n=1,\qquad 1+a_{j-1}<ca_j.
\]
They are feasible exactly when
$c^n>1+c+\cdots+c^{n-1}$, so their thresholds are the $\beta_n$
of Section~\ref{sec:two}, with $\beta_n\uparrow2$.
With mass $p$ at zero and equal masses $(1-p)/n$ at the positive
atoms, counting guaranteed successes and replacing zero recursively
gives
\[
 \frac{p(1-p)(2-p)
       -p(1-p)^3\left(\frac{3}{2n}+\frac{1}{2n^2}\right)}
      {1-p^4}.
\]
As $n\to\infty$ this approaches the value of \citet{BF} after
optimizing $p$. The plotted finite witnesses are checked by direct
four-draw enumeration; numerical mass optimization does not assert
their optimality.

Section~\ref{sec:four-exact} gives the exact endpoint
$S_3(2)=C_{\mathrm{BF}}=0.4006953723293\ldots$.
The blue bounds below $2$, described below, include
$S_3(19/10)<0.381302383<C_{\mathrm{BF}}$.
Neither these bounds nor the finite ladders resolve whether the
value $C_{\mathrm{BF}}$ holds on an interval immediately below $2$.
The figure therefore draws a green point at $2$, rather than an
exact interval ending there.

\subsubsection{Above two: A plateau, a window, and a jump}
For $2<c<3$, recursive zero replacement of a single positive atom
gives Yu's lower bound \cite{Yu}
\[
 L^*=\sup_{w\in(0,1)}\frac{w-w^4}{1-(1-w)^4}
     =0.482781746628666\ldots.
\]
Writing $z=1-w$, the function to maximize is
\[
 L(z)=\frac{3z-3z^2+z^3}{1+z+z^2+z^3}.
\]
Its derivative has numerator $4z^4-4z^3-3z^2-6z+3$,
which decreases strictly from $3$ to $-6$ on $[0,1]$.
Its unique zero $z_*=0.3963799474937\ldots$ therefore gives
$L^*=L(z_*)$.
Put $\gamma=(1+\sqrt{13})/2$. The first exact interval above $2$ is
\begin{equation}\label{eq:s3-first-plateau}
 S_3(c)=L^*\qquad(2<c\le\gamma).
\end{equation}

\begin{proof}
For the lower bound, take
$\Prob[X=\rho^k]=(1-z_*)z_*^k$, $k\ge0$, with
$0<\rho<\min\{(c-2)/3,1/c\}$.
The success probability satisfies
$P=(w-w^4)+z_*^4P$, with $w=1-z_*$, and hence equals $L^*$.
At the endpoint $\gamma$, record the pair cutoffs
\[
 \{\gamma/2,\gamma-1,\gamma/(\gamma-1),\gamma\}
\]
and all three potentially winning target positions.
The ordinary four-point certificate has exact score at most $L^*$
on every pattern and positive semidefinite matrices over the quartic
field containing $L^*$. Thus $S_3(\gamma)\le L^*$.
Monotonicity proves the interval.
Appendix~\ref{app:s3-certificates} describes the complete cover and
the exact matrix checks.
\end{proof}

Consequently $\lim_{c\to2^+}S_3(c)=L^*$, answering Yu's question,
and $S_3$ jumps from $C_{\mathrm{BF}}$ to $L^*$ immediately after $2$.
The next change in the construction comes from two positive atoms
$r,1$. To make $(r,r,r)$ win against $1$ and $(0,r,1)$ win
against $r$, we need
\[
 3r<c,\qquad 1+r<cr,
 \quad\text{or equivalently}\quad
 \frac1{c-1}<r<\frac c3.
\]
This window opens exactly when $c^2-c-3>0$, that is, $c>\gamma$.
The window opening provides the left threshold of the next exact
interval; its right endpoint will again be certified separately.

\begin{proposition}\label{prop:three-above}
Let $\nu=(1-w)\delta_0+wp\,\delta_r+w(1-p)\,\delta_1$.
\begin{enumerate}
  \item If $\gamma<c<3$, $1/(c-1)<r<c/3$, $p=7/8$ and $w=5/8$,
  recursive zero replacement gives
  $S_3(c)\ge146857/299008>0.4911473$.
  \item If $1+\sqrt2<c<3$, $1/(c-1)<r<(c-1)/2$,
  $p=29/40$ and $w=13/20$, it gives
  $S_3(c)\ge7808829037/15517440000>0.5032292>1/2$.
\end{enumerate}
Both bounds extend to $c\ge3$ by monotonicity.
\end{proposition}

\begin{proof}
Put $z=1-w$, $a=wp$ and $b=w(1-p)$.
In the first window, the guaranteed winning configurations have
total probability
\[
\begin{aligned}
 N_1={}&a\bigl[z^3+3z^2(a+b)+3z(a^2+2ab)\bigr]\\
       &+b\bigl[z^3+3z^2(a+b)+3z(a+b)^2+a^3\bigr].
\end{aligned}
\]
Against target $r$, these are the triples with at least two zeros,
and those with one zero except $(0,1,1)$. Against target $1$,
every triple with a zero and $(r,r,r)$ wins.
The second window also makes $(r,r,1)$ win against $1$, giving
$N_2=N_1+3a^2b^2$.
The four-zero event repeats the construction, so
$S_3(c)\ge N_i/(1-z^4)$, by \cite[Proposition~3.3]{BF}, which applies to any number
of draws, with scales separated enough that lower scales do not change the counted comparisons. Substitution gives the stated fractions.
The second window is nonempty exactly when $c^2-2c-1>0$,
so its threshold is the silver ratio $1+\sqrt2$.
\end{proof}

For $z,a,b>0$ with $z+a+b=1$, the first numerator simplifies to
\[
 N(z,a,b)=(a+b)-(a+b)^4+ab(a^2-3zb).
\]
Optimizing the same two-positive-atom family gives
\begin{equation}\label{eq:s3-second-value}
 L_2=\max_{\substack{z,a,b>0\\z+a+b=1}}
       \frac{N(z,a,b)}{1-z^4}
     =0.491156081813797\ldots.
\end{equation}
Its masses are $z=0.3767059973499\ldots$,
$a=0.5439163904458\ldots$ and $b=0.0793776122041\ldots$,
represented exactly in a degree-$14$ field.
An ordinary four-point certificate at $47/20=2.35$, adding the
cutoff $3/c$ to the four basic cutoffs, proves
$S_3(47/20)\le L_2$. It also proves global optimality of the stated
masses within the construction.
Together with the uniform lower bound for $c>\gamma$, this gives
\begin{equation}\label{eq:s3-second-plateau}
 S_3(c)=L_2\qquad(\gamma<c\le47/20).
\end{equation}
Thus the jump immediately after $\gamma$ is $L_2-L^*>0.008374335$.
The endpoint $\gamma$ is sharp for the first plateau; $47/20$ is
the proved extent of the second one, whose maximal right endpoint
remains unknown.

\subsubsection{Finite ladders approaching the silver ratio}
The next orange steps in the lower panel come from inserting more
positive atoms. Before normalizing the top atom to $1$, use
\[
 1=a_0<a_1<\cdots<a_N<c-1,\qquad
 2+a_j<ca_{j+1}\quad(0\le j<N).
\]
The second comparison makes $(1,1,a_j)$ win against $a_{j+1}$.
The upper endpoint ensures that $(0,1,a_j)$ wins against $1$.
The limiting positions obey
\[
 r_0=1,\qquad r_{j+1}=\frac{2+r_j}{c},\qquad
 r_N=\frac{3+2\sum_{k=1}^{N-1}c^k}{c^N}.
\]
Every strict ladder has $a_N>r_N$, and if $r_N<c-1$,
small upward perturbations give a strict ladder.
Its feasibility threshold $\theta_N$ therefore satisfies
\[
 \theta_N^{N+1}-\theta_N^N
   =3+2\sum_{k=1}^{N-1}\theta_N^k.
\]
Here $\theta_1=\gamma$, and $\theta_N$ increases to $1+\sqrt2$:
The recurrence has fixed point $2/(c-1)$, which meets $c-1$
at $(c-1)^2=2$. These are construction thresholds; we do not
claim that every $\theta_N$ is a discontinuity of $S_3$.

Put mass $p$ at zero, $q$ at $1$, and $r/N$ at each of the
$N$ tail atoms, with $p+q+r=1$.
Counting guaranteed successes gives
\[
\begin{aligned}
 S_3(c)&\ge\frac{F_N(p,q,r)}{1-p^4},\\
 F_N&=p^3(q+r)+3p^2(q+r)^2
      +3p(q^3+3q^2r+2qr^2+r^3)\\
    &\qquad+q^3r+\tfrac32q^2r^2(1-\tfrac1N).
\end{aligned}
\]
The terms group configurations by the number of zero summands;
the last two terms count $(1,1,1)$ and $(1,1,a_j)$ against a
larger tail target. The denominator repeats the four-zero event
on a smaller scale. Fixed rational masses and a sufficiently small
scale ratio turn this count into a verified lower witness.

Figure~\ref{fig:three-against-one} uses ten saved ladder lengths.
For each $N$, an exact rational calculation brackets $\theta_N$
to width $10^{-9}$, and the witness is placed at the feasible upper
end of that bracket. It therefore starts just after its threshold,
rather than at the coefficient where numerical mass search was run.
For example, the two- and three-tail-atom witnesses start at
$2.374423764$ and $2.398726488$, with lower values exceeding
$0.492838853$ and $0.493537310$. The ten-tail-atom witness exceeds
$0.494663204$ from $2.414182765$ onward.
These starts and probabilities are checked in
\path{lower_thresholds.json}.

At the silver ratio itself, finite ladders of every length are
feasible below $\sqrt2$. Taking $N\to\infty$ in the same polynomial
and choosing the algebraic masses specified in
Appendix~\ref{app:s3-lower} gives
\[
 S_3(1+\sqrt2)\ge L_3=0.4952063382\ldots,
\]
where $L_3$ is algebraic of degree $14$.
This is the orange endpoint value in the figure; exact enumeration
also checked finite ladders of lengths $1,3,4,5$ over
$\mathbb Q(\sqrt2)$.

\subsubsection{Certifying the blue curve}
On four sorted points, record each potentially winning target
position and the pair comparisons $x_j<q x_i$.
For $c<3$, the smallest point cannot win. For $c>2$,
the second-smallest can win, so all three other targets must be
retained. The basic cutoffs above $2$ are
\[
 \{c/2,c-1,c/(c-1),c\}.
\]
We add $3/c$ for the second exact interval and the bound at $2.36$.
At $2.4$, the limiting ladder positions provide three further
cutoffs, $5/4$, $65/48$ and $805/576$.
The certificates bound the objective plus rooted flag squares on
every retained pattern; excluded prefixes carry exact Gordan
certificates. Taking expectations then gives an upper bound
for every law, as in Section~\ref{sec:pair-colors-general}.

Below $2$, four-point certificates give
\[
\begin{aligned}
 S_3(3/2)&<0.3221953,&
 S_3(8/5)&<0.3371723,&
 S_3(7/4)&<0.3532920,\\
 S_3(37/20)&<0.3738948,&
 S_3(19/10)&<0.381302383.
\end{aligned}
\]
The first and third bounds leave gaps of less than $0.000826$
and $0.002392$ from the verified lower witnesses.
The cutoff sets vary across this range.
At $3/2$, $8/5$ and $7/4$, they are drawn from
$\{2/c,c,1/(c-1),3/c,c/(c-1)\}$.
The certificate at $37/20$ instead uses
$\{2/c,7/6,3/2,c,2,c/(c-1)\}$.
Above the second exact interval, we have
\begin{equation}\label{eq:s3-colored}
\begin{aligned}
 S_3(59/25)&\le
   \frac{14739078469921721}{30000000000000000}<0.491302616,\\
 S_3(12/5)&\le
   \frac{7403289897571463}{15000000000000000}<0.493552660,\\
 S_3(121/50)&\le
   \frac{10158973613564103}{20000000000000000}<0.507948681.
\end{aligned}
\end{equation}
The gap at $2.4$ to its rational ladder witness is less than
$1.54\cdot10^{-5}$. Here recording the ladder ratios improves
the certificate while keeping four sample points.

At the exact silver-ratio endpoint, the stronger upper bound
still comes from seven points:
\begin{equation}\label{eq:s3-seven}
 S_3(1+\sqrt2)\le
 \frac{96556198460261613624314242853}{194975556187447393296612065280}
 <0.4952221.
\end{equation}
The pattern cover is checked over $\mathbb Q(\sqrt2)$ and the
PSD matrices are rational Gram factors.
Appendix~\ref{app:s3-certificates} records the certificate sizes
and the two independent verification routes.
The answer to Problem~5.3 needs less: a four-point certificate of the
same form as the one at $\gamma$, with cutoffs $c/2$, $c-1$, $c/(c-1)$ and $c$,
already gives $S_3(1+\sqrt2)<0.4988589<\tfrac12$
(\path{anc/proofs/s3/silver_ratio_four_point/}).

All four-point data are in
\path{anc/proofs/s3/}.
The standalone scripts \path{verify_below_two.py},
\path{verify_endpoint.py} and \path{verify_second_plateau.py}
check the three nontrivial exact intervals.
The script \path{verify_all.py} checks the saved rational upper
certificates, lower witnesses, threshold brackets and uniform
window bounds. Each lower witness contains a rational support,
masses and a scale ratio preserving its counted successes.
These exact checks certify the curves; numerical solver residuals
and numerical mass optimization are not proof inputs.

The plotting script \path{figures/make_three_against_one.py} takes
the smallest available upper bound with endpoint at least $c$ and
the largest lower bound whose witness is valid at $c$.
It discards bounds dominated by monotonicity and draws each exact
interval once. This explains both the horizontal bound segments
and the shaded bands. In particular, the seven-point bound applies
up to and including the silver ratio, whereas the certificate at
$2.42$ supplies an upper bound after it only through $2.42$.

\subsubsection{The silver-ratio endpoint and fairness}
The bounds at the silver ratio give
\[
 0.4952063<S_3(1+\sqrt2)<0.4952221<\tfrac12.
\]
Immediately to its right, the second window in
Proposition~\ref{prop:three-above} opens.
For every $1+\sqrt2<c<3$, choose
$1/(c-1)<r<(c-1)/2$ and use the fixed rational masses
\[
 (z,a,b)=10^{-12}(339588657305,\;481077629994,\;179333712701)
\]
at $0,r,1$. The same success polynomial gives
\begin{equation}\label{eq:s3-uniform-lower}
 \frac{N_1+3a^2b^2}{1-z^4}
 =0.503419372374\ldots>0.503419372.
\end{equation}
The masses sum to one and work throughout the window, with
$r$ and the scale separation depending on $c$.
Monotonicity extends the bound to $c\ge3$.
The exact polynomial check and a witness at $121/50$ are saved
in \path{lower_uniform.json}.
Thus the orange curve immediately after $1+\sqrt2$ lies above
$1/2$, while the blue endpoint value lies below it.

\begin{proof}[Proof of Theorem~\ref{thm:fair3}]
The constructions and matching endpoint certificates above prove
the exact intervals and the jumps at $\sqrt2$ and $\gamma$.
The ladder lower bound and \eqref{eq:s3-seven} give the stated
silver-ratio band.
By monotonicity, \eqref{eq:s3-seven} implies $S_3(c)<1/2$
for every $c\le1+\sqrt2$.
The uniform construction \eqref{eq:s3-uniform-lower} gives
$S_3(c)>1/2$ for every $c>1+\sqrt2$.
Hence no positive coefficient makes the house's best chance
exactly even.
\end{proof}

\addtocontents{toc}{\protect\setcounter{tocdepth}{1}}

%======================================================================
\section{Open problems}\label{sec:open}
%======================================================================

\subsection{More discontinuity points and exact values/plateaus}

A first set of open problems is to make the staircase picture
in Table~\ref{tab:steps-bounds} and Figures~\ref{fig:staircase} and~\ref{fig:three-against-one}
more precise.

\subsubsection{Certify more exact values and plateaus of $S_2$ and $S_3$}
Natural candidates are the remaining thresholds of Table~\ref{tab:steps-thresholds} where a discontinuity is proved but the exact value is not yet certified.

\subsubsection{Certify more discontinuity points of $S_2$ and $S_3$}
Natural candidates here are the thresholds $c_n$ for $n>6$,
$\tau_k$ for $k>3$, and $\alpha_m$ for $m>4$
where the sequences $(c_n)$, $(\tau_k)$ and $(\alpha_m)$
extend the definitions in Table~\ref{tab:steps-thresholds}
in the natural way.

\subsubsection{Find an infinite number of discontinuities of $S_2$ and $S_3$}
Here, natural candidates are the sequences $(c_n)$, $(\tau_k)$ and $(\alpha_m)$. While the previous two problems can be solved by finite computations, this one requires a more conceptual argument.

\subsubsection{Determine all discontinuity points near an
    accumulation point
} For instance, for $S_2$ near $\varphi$ and for $S_3$ near $1+\sqrt2$.
Knowing that $c_4 < \tau_1 < c_5 < \tau_2 < c_6 < \tau_3$
are discontinuity points does not rule out other discontinuity points
between them.

\subsection{On the continuous part of $S_3$ and/or $S_2$}
Both functions are nondecreasing and left-continuous, so they have at most countably many discontinuities. Nondecreasing functions can be decomposed
into the sum of a continuous nondecreasing part and a jump part.
Our results so far suggest that the continuous part is 0,
and we prove this in $(0,c_3)$ and $(\varphi, \alpha_3)$ for $S_2$ and in $(0, \sqrt2)$ and $(2, 47/20)$ for $S_3$
(see Figures~\ref{fig:staircase} and~\ref{fig:three-against-one}).

\subsection{Relation to the one-scale problem and $V(\beta)$}
Is $S_2(c)=V(\beta(c))$, where $\beta(c)$ is the one-scale problem $\sup\{\Prob_\rho[Y_1+Y_2<cY_3]:\ \mathrm{supp}\,\rho\subseteq[a,1],\ ca>1\}$? Is $S_2$ a step function, with steps at algebraic numbers where new windows open?
The value $(2\sqrt7-4)/3$ joins a family of quadratic-irrational extremal values with self-similar extremizers; for example, the packing density of the permutation $132$ is $2\sqrt3-3$, a result of Stromquist and of Galvin and Kleitman (see \cite{PS10}).
\subsection{Beyond linear inequalities}
The method applies to any event cut out by strict linear inequalities, for instance the other members of the family \cite[(3.4)]{BF}. For events defined by polynomial inequalities, such as the obtuse random triangles of Lewis Carroll's pillow problem discussed in \cite[Section~5]{BF} and \cite{DR25}, the patterns are still sign vectors, but certifying which of them are realizable requires semialgebraic rather than linear arguments.

\subsection{Self-similarity}
Some optimal laws $\mu$ used in lower bound constructions
are self-similar. One aspect of this is the use
of Proposition~3.3 of \citet{BF} which repeats a starting measure
$\nu$ at infinitely many smaller scales to obtain a measure
$\mu$ satisfying \eqref{eq:boost}.
Another aspect of self-similarity is
Lemma~\ref{lem:subst4} which substitutes the mass of $\mu$ below $b$
by a copy of $\mu$ scaled by a small $\varepsilon>0$.

In extremal combinatorics, exact results with self-similar extremizers were also obtained by combining flag-algebra certificates with stability \cite{PST19} and a recursion over the levels of an iterated blow-up \cite{BHLP16}.
But deeper connections between these known results in extremal combinatorics
and the self-similarity of some constructions in the present paper
are unclear at this point.

\addtocontents{toc}{\protect\setcounter{tocdepth}{-1}}
\section{Code and data, acknowledgements}
\subsection*{Back at the casino}
%We started with a casino that can make the three-cup game unfavourable for you. 
In the four-roll game of the introduction, the optimal probability of $X_1+X_2+X_3<2X_4$ is exactly the value $C_{\mathrm{BF}}=0.4006953\ldots$ conjectured by \citet{BF}. In the two-against-one game, the house can choose both the dice and the coefficient $c$, yet it cannot make its best chance exactly even: up to and including the golden ratio, its best chance stays below $1/2$, and for every larger coefficient it exceeds $1/2$. When you roll three times against the house's one, the silver ratio $1+\sqrt2$ plays the role of the golden ratio, and again no coefficient makes the house's best chance exactly even; this answers Problem~5.3 of \cite{BF}.

\subsection*{Code and data}
The certificates, exact verifiers, lower constructions and Lean proof are
available in the public repository
\url{https://github.com/antalez/flag-algebras-iid}, under \path{anc/},
and in the arXiv ancillary files with the same structure but without
the \path{anc/} prefix. Table~\ref{tab:ancillary-files} lists the directories
using the arXiv paths; prepend \path{anc/} for their repository locations.

Result subdirectories are ordered by their lower values, truncated to
nine decimal places, followed by the interval or endpoint name.
Table~\ref{tab:ancillary-s2} lists the subdirectories of \path{proofs/s2/}
and their intervals.
For example, \path{0.449989622__c1_c2/} contains the proof for
$(c_1,c_2]$. The ancillary \path{README.md} gives verification instructions,
\path{results.json} indexes the results and their files, and
\path{verify.py} lists and runs the relocated Python certificate checks.
The six-point certificates at $\varphi$, $\beta_3'$ and $\beta_3$ are in
\path{proofs/s2/six_point_endpoints/}.

\begin{table}[htbp]
\centering
\small
\begin{tabular}{@{}p{0.46\linewidth}p{0.50\linewidth}@{}}
\toprule
Directory & Contents\\
\midrule
\path{proofs/common/}
  & Shared exact arithmetic and its tests.\\
\path{proofs/s2/}
  & Result subdirectories for the $S_2$ bounds and constructions in
    Tables~\ref{tab:steps-thresholds}--\ref{tab:steps-bounds}.\\
\path{proofs/s2/common/}
  & Colored-flag and near-optimal-law verifiers.\\
\path{proofs/s3/}
  & Result subdirectories for the exact $S_3$ intervals and the
    four-roll equality.\\
\path{proofs/s3/0.400695372__c2/}
  & The exact certificate for $C=C_{\mathrm{BF}}$ and its two
    independent Python verifiers.\\
\path{proofs/s3/common/}
  & Colored-flag verifiers and lower-witness checks.\\
\path{proofs/s3/rational/}
  & Upper certificates at rational coefficients for Figure~\ref{fig:three-against-one}.\\
\path{proofs/s3/lower_constructions/}
  & Rational lower witnesses for Figure~\ref{fig:three-against-one}.\\
\path{proofs/s3/silver_ratio_seven_point/}
  & The seven-point silver-ratio verifier and compressed exclusion witnesses;
    the large Gram-factor input is available in the repository.\\
\path{lean/}
  & The Lean proof of Theorem~\ref{thm:main}, including its shared
    dependencies and earlier proof routes.\\
\bottomrule
\end{tabular}
\caption{Directory structure of the ancillary files. Repository paths
have the additional prefix \texttt{anc/}.}
\label{tab:ancillary-files}
\end{table}

\begin{table}[htbp]
\centering
\small
\begin{tabular}{@{}p{0.46\linewidth}p{0.19\linewidth}p{0.29\linewidth}@{}}
\toprule
Subdirectory of \path{proofs/s2/} & Interval & Contents\\
\midrule
\path{0.430500874__1_c1/}
  & $(1,c_1]$ & Exact plateau\\
\path{0.449989622__c1_c2/}
  & $(c_1,c_2]$ & Exact plateau\\
\path{0.457361518__c2_c3/}
  & $(c_2,c_3]$ & Exact plateau\\
\path{0.460603680__c3_tau0/}
  & $(c_3,\tau_0]$ & Certified band\\
\path{0.461389881__tau0_c4/}
  & $(\tau_0,c_4]$ & Certified band\\
\path{0.462747996__c4_tau1/}
  & $(c_4,\tau_1]$ & Certified band\\
\path{0.463221002__tau1_c5/}
  & $(\tau_1,c_5]$ & Certified band\\
\path{0.463950398__c5_tau2/}
  & $(c_5,\tau_2]$ & Certified band\\
\path{0.464266822__tau2_c6/}
  & $(\tau_2,c_6]$ & Certified band\\
\path{0.464719509__c6_tau3/}
  & $(c_6,\tau_3]$ & Certified band\\
\path{0.464946735__tau3_phi/}
  & $(\tau_3,\varphi]$ & Lower construction\\
\path{0.503098408__phi_alpha1/}
  & $(\varphi,\alpha_1]$ & Exact plateau\\
\path{0.507073883__alpha1_alpha2/}
  & $(\alpha_1,\alpha_2]$ & Exact plateau\\
\path{0.508460361__alpha2_alpha3/}
  & $(\alpha_2,\alpha_3]$ & Exact plateau\\
\path{0.509054132__alpha3_alpha4/}
  & $(\alpha_3,\alpha_4]$ & Certified band\\
\path{0.509376431__alpha4_beta3prime/}
  & $(\alpha_4,\beta_3']$ & Lower construction\\
\path{common/}
  & All intervals & Shared verifiers\\
\bottomrule
\end{tabular}
\caption{Subdirectories of \texttt{proofs/s2/}, ordered by the nine-digit
lower-value prefix. A certified band has a positive gap between the
lower and upper bounds. For the two directories containing only a lower
construction, the upper certificates are in \texttt{proofs/s2/six\_point\_endpoints/},
as are the files for $(\beta_3',\beta_3]$.
The elementary plateau $(0,1]$ requires no certificate files.}
\label{tab:ancillary-s2}
\end{table}

\subsection*{Acknowledgements}
We thank Tobias Fritz for posing the question on MathOverflow \cite{MO474916}, for his encouragement from the first draft on, and for suggesting that we present the approach as a general method. He also asked whether the certificates reflect the separation of scales in the near-optimal constructions. Section~\ref{sec:four-exact} answers this question for the four-point certificate.

\appendix
\section{Completeness of Table~\ref{tab:example}}\label{app:gordan}

This appendix shows how the fourteen patterns of Table~\ref{tab:example} were found and why there are no others. Here $c=6/5$. A depth-first search decides the eight relations of Section~\ref{sec:example} one at a time, in the order listed there, and abandons a partial sign vector as soon as it is impossible. Each line below gives an abandoned prefix and nonnegative weights on its relations, written $y\,(t)$ for weight $y$ on relation $t$. To check a line, multiply relation $t$ by its sign ($+$ for ``holds'', $-$ for ``fails'') so that it reads ``form $<0$''. The weighted sum of these forms, written in the gaps $d_a>0$, then has all coefficients $\ge0$, so it cannot be negative, and the prefix is impossible. For example, the prefix \verb|+-| says $x_1+x_2<cx_3$ and $x_1+x_2>cx_4$. Half of the first plus half of the second gives $\tfrac c2(x_4-x_3)<0$, which is impossible. The $14$ sign vectors that survive are those of Table~\ref{tab:example}.

\begin{center}
\small
\begin{tabular}{ll@{\qquad}ll}
\toprule
prefix & certificate & prefix & certificate\\
\midrule
\verb|+-| & $\tfrac{1}{2}\,(1) + \tfrac{1}{2}\,(2)$ & \verb|++++-++| & $\tfrac{1}{6}\,(3) + \tfrac{5}{11}\,(6) + \tfrac{25}{66}\,(7)$\\
\verb|-++| & $\tfrac{1}{2}\,(1) + \tfrac{1}{2}\,(3)$ & \verb|++----+| & $\tfrac{1}{2}\,(4) + \tfrac{1}{2}\,(7)$\\
\verb|--+| & $\tfrac{1}{2}\,(1) + \tfrac{1}{2}\,(3)$ & \verb|-+----+| & $\tfrac{1}{2}\,(4) + \tfrac{1}{2}\,(7)$\\
\verb|+++-| & $\tfrac{3}{11}\,(1) + \tfrac{5}{22}\,(3) + \tfrac{1}{2}\,(4)$ & \verb|------+| & $\tfrac{1}{2}\,(2) + \tfrac{1}{2}\,(7)$\\
\verb|---+| & $\tfrac{1}{2}\,(2) + \tfrac{1}{2}\,(4)$ & \verb|++++++-+| & $\tfrac{1}{6}\,(5) + \tfrac{5}{6}\,(8)$\\
\verb|++-++| & $\tfrac{1}{2}\,(3) + \tfrac{1}{2}\,(5)$ & \verb|++++-+-+| & $\tfrac{209}{723}\,(2) + \tfrac{275}{723}\,(3) + \tfrac{59}{723}\,(7) + \tfrac{60}{241}\,(8)$\\
\verb|++--+| & $\tfrac{1}{2}\,(3) + \tfrac{1}{2}\,(5)$ & \verb|++++--++| & $\tfrac{1}{6}\,(3) + \tfrac{5}{6}\,(8)$\\
\verb|-+-++| & $\tfrac{1}{2}\,(1) + \tfrac{1}{2}\,(5)$ & \verb|++++---+| & $\tfrac{209}{723}\,(2) + \tfrac{275}{723}\,(3) + \tfrac{59}{723}\,(7) + \tfrac{60}{241}\,(8)$\\
\verb|-+--+| & $\tfrac{1}{2}\,(1) + \tfrac{1}{2}\,(5)$ & \verb|++-+-+-+| & $\tfrac{1}{2}\,(7) + \tfrac{1}{2}\,(8)$\\
\verb|----+| & $\tfrac{1}{2}\,(1) + \tfrac{1}{2}\,(5)$ & \verb|++-+--++| & $\tfrac{1}{2}\,(6) + \tfrac{1}{2}\,(8)$\\
\verb|+++++-| & $\tfrac{1}{2}\,(5) + \tfrac{1}{2}\,(6)$ & \verb|++-+---+| & $\tfrac{1}{4}\,(4) + \tfrac{1}{4}\,(6) + \tfrac{1}{4}\,(7) + \tfrac{1}{4}\,(8)$\\
\verb|++---+| & $\tfrac{1}{2}\,(4) + \tfrac{1}{2}\,(6)$ & \verb|++-----+| & $\tfrac{1}{2}\,(6) + \tfrac{1}{2}\,(8)$\\
\verb|-+-+-+| & $\tfrac{1}{2}\,(1) + \tfrac{5}{22}\,(2) + \tfrac{3}{11}\,(6)$ & \verb|-+-+--++| & $\tfrac{1}{7}\,(1) + \tfrac{6}{7}\,(8)$\\
\verb|-+---+| & $\tfrac{1}{2}\,(1) + \tfrac{1}{2}\,(6)$ & \verb|-+-+---+| & $\tfrac{11}{34}\,(1) + \tfrac{11}{34}\,(2) + \tfrac{3}{17}\,(7) + \tfrac{3}{17}\,(8)$\\
\verb|-----+| & $\tfrac{1}{2}\,(2) + \tfrac{1}{2}\,(6)$ & \verb|-+-----+| & $\tfrac{11}{34}\,(1) + \tfrac{11}{34}\,(2) + \tfrac{3}{17}\,(4) + \tfrac{3}{17}\,(8)$\\
\verb|+++++++| & $\tfrac{5}{9}\,(5) + \tfrac{4}{9}\,(7)$ & \verb|-------+| & $\tfrac{1}{7}\,(2) + \tfrac{6}{7}\,(8)$\\
\bottomrule
\end{tabular}
\end{center}

% Completeness at c = phi: Gordan certificates in the format of Appendix A.
% Generated and checked in exact arithmetic over Q(sqrt 5); relations numbered as in Section 2.4.
\section{Four-point pattern completeness at the golden ratio}\label{app:gordan-phi}
Use the conventions of Appendix~\ref{app:gordan}, now with
$c=\varphi=(1+\sqrt5)/2$. The certificates below exclude all but
$16$ sign vectors, in the order of the eight relations of
Section~\ref{sec:example}:
\[
\begin{gathered}
\texttt{+++++++-},\quad \texttt{++++++--},\quad
\texttt{++++-++-},\quad \texttt{++++-+--},\\
\texttt{++++--+-},\quad \texttt{++++----},\quad
\texttt{++-+-+++},\quad \texttt{++-+-++-},\\
\texttt{++-+-+--},\quad \texttt{++-+--+-},\quad
\texttt{++-+----},\quad \texttt{++------},\\
\texttt{-+-+--+-},\quad \texttt{-+-+----},\quad
\texttt{-+------},\quad \texttt{--------}.
\end{gathered}
\]
\begin{center}
\small
\begin{tabular}{ll@{\qquad}ll}
\toprule
prefix & certificate & prefix & certificate\\
\midrule
\verb|+-| & $\tfrac{1}{2}\,(1) + \tfrac{1}{2}\,(2)$ & \verb|++-----+| & $\tfrac{1}{2}\,(4) + \tfrac{1}{2}\,(8)$\\
\verb|+++-| & $\tfrac{1}{2}\,(3) + \tfrac{1}{2}\,(4)$ & \verb|-++| & $\tfrac{1}{2}\,(1) + \tfrac{1}{2}\,(3)$\\
\verb|+++++-| & $\tfrac{1}{2}\,(5) + \tfrac{1}{2}\,(6)$ & \verb|-+-++| & $\tfrac{1}{2}\,(1) + \tfrac{1}{2}\,(5)$\\
\verb|++++++++| & $\tfrac{3-\sqrt5}{2}\,(3) + \tfrac{-1+\sqrt5}{2}\,(8)$ & \verb|-+-+-+| & $\tfrac{1}{2}\,(1) + \tfrac{1}{2}\,(6)$\\
\verb|++++++-+| & $\tfrac{1}{2}\,(7) + \tfrac{1}{2}\,(8)$ & \verb|-+-+--++| & $\tfrac{1}{2}\,(1) + \tfrac{1}{2}\,(8)$\\
\verb|++++-+++| & $\tfrac{3-\sqrt5}{2}\,(3) + \tfrac{-1+\sqrt5}{2}\,(8)$ & \verb|-+-+---+| & $\tfrac{1}{2}\,(1) + \tfrac{1}{2}\,(8)$\\
\verb|++++-+-+| & $\tfrac{1}{2}\,(7) + \tfrac{1}{2}\,(8)$ & \verb|-+--+| & $\tfrac{1}{2}\,(1) + \tfrac{1}{2}\,(5)$\\
\verb|++++--++| & $\tfrac{1}{2}\,(6) + \tfrac{1}{2}\,(8)$ & \verb|-+---+| & $\tfrac{1}{2}\,(1) + \tfrac{1}{2}\,(6)$\\
\verb|++++---+| & $\tfrac{1}{2}\,(6) + \tfrac{1}{2}\,(8)$ & \verb|-+----+| & $\tfrac{1}{2}\,(4) + \tfrac{1}{2}\,(7)$\\
\verb|++-++| & $\tfrac{1}{2}\,(3) + \tfrac{1}{2}\,(5)$ & \verb|-+-----+| & $\tfrac{1}{2}\,(1) + \tfrac{1}{2}\,(8)$\\
\verb|++-+-+-+| & $\tfrac{1}{2}\,(7) + \tfrac{1}{2}\,(8)$ & \verb|--+| & $\tfrac{1}{2}\,(1) + \tfrac{1}{2}\,(3)$\\
\verb|++-+--++| & $\tfrac{1}{2}\,(6) + \tfrac{1}{2}\,(8)$ & \verb|---+| & $\tfrac{1}{2}\,(2) + \tfrac{1}{2}\,(4)$\\
\verb|++-+---+| & $\tfrac{1}{2}\,(6) + \tfrac{1}{2}\,(8)$ & \verb|----+| & $\tfrac{1}{2}\,(1) + \tfrac{1}{2}\,(5)$\\
\verb|++--+| & $\tfrac{1}{2}\,(3) + \tfrac{1}{2}\,(5)$ & \verb|-----+| & $\tfrac{1}{2}\,(1) + \tfrac{1}{2}\,(6)$\\
\verb|++---+| & $\tfrac{1}{2}\,(4) + \tfrac{1}{2}\,(6)$ & \verb|------+| & $\tfrac{1}{2}\,(2) + \tfrac{1}{2}\,(7)$\\
\verb|++----+| & $\tfrac{1}{2}\,(4) + \tfrac{1}{2}\,(7)$ & \verb|-------+| & $\tfrac{1}{2}\,(1) + \tfrac{1}{2}\,(8)$\\
\bottomrule
\end{tabular}
\end{center}

\section{Certificate details for three against one}\label{app:s3-certificates}
This appendix records the pattern covers, matrix dimensions and
exact arithmetic used by the upper certificates in
Section~\ref{sec:three}. The main argument needs the complete
pattern cover, PSD positivity, and the score inequality on every
retained pattern.

\subsection{The endpoint below two}
At $c=\sqrt2$, the four-point cover retains $24$ patterns and has
$46$ exact Gordan exclusions over $\mathbb Q(\sqrt2)$.
The two $(2,3)$ blocks have flag dimensions $8,7$ and are stored as
$Q_\sigma=W_\sigma Y_\sigma W_\sigma^\top$, with positive definite
$Y_\sigma$ of dimensions $2,3$ over the cubic field $\mathbb Q(u)$,
where $u=1/z$ and $z$ is the optimizer in \eqref{eq:s3-below-two}.
The matrices are
\[
 Y_0=\begin{pmatrix}D&1/4\\1/4&1\end{pmatrix},
 \qquad Y_1=\operatorname{diag}(1/2,1/2,L_0),
 \qquad D=-\frac{13}{4}+\frac{3u}{4}+\frac{3u^2}{2}.
\]
The root satisfies $u>7/5$, so $D>1/16$ and both matrices are
positive definite. All scores are at most $L_0$, with equality on
nine patterns. The script \path{verify_below_two.py} checks the
cover, positive definiteness and every exact score.

\subsection{The two exact intervals above two}
At $\gamma=(1+\sqrt{13})/2$, the four-point cover retains $493$
patterns; $2{,}535$ exact Gordan certificates over
$\mathbb Q(\sqrt{13})$ check every excluded prefix.
There is one $(0,2)$ block and five $(2,3)$ blocks,
with flag dimensions $5,27,23,27,19,21$.
Their matrices are $Q_\sigma=W_\sigma Y_\sigma W_\sigma^\top$,
where the $Y_\sigma$ have dimensions $2,14,22,26,15,16$ and are
positive definite over $\mathbb Q(u)$.
Here $u=1/z_*$ is the root in $(2,3)$ of
$3u^4-6u^3-3u^2-4u+4=0$.
Exact LDL elimination verifies their positivity.
Every pattern score is at most $L^*$, with equality on $62$ patterns.
The standard-library script \path{verify_endpoint.py} rebuilds
the cover, flag products and scores independently of the numerical
search.

At $47/20$, the second-plateau cover retains $931$ patterns and has
$5{,}645$ rational Gordan exclusions.
There is one $(0,2)$ block and six $(2,3)$ blocks, with original
dimensions $6$, $33$, $32$, $35$, $25$, $23$, $29$ and reduced
dimensions $2$, $16$, $31$, $34$, $19$, $19$, $20$.
The matrices again have the form
$Q_\sigma=W_\sigma Y_\sigma W_\sigma^\top$.
Rational LDL elimination together with exact algebraic entry
enclosures verifies that each $Y_\sigma$ is positive definite.
Every score is at most the exact degree-$14$ value $L_2$ in
\eqref{eq:s3-second-value}, with equality on $120$ patterns.
The standard-library script \path{verify_second_plateau.py}
checks all products and scores independently of the numerical search.
These scripts and their proof inputs are in
\path{anc/proofs/s3/}.

\subsection{Rational four-point bounds}
The covers at $59/25$ and $121/50$ retain $934$ and $571$ patterns,
respectively. Their certificates have eight and seven PSD blocks,
of dimension at most $35$. At $12/5$, the ladder-refined cover retains
$2{,}625$ patterns and has $18$ PSD blocks, of dimension at most $87$.
The exact rational Gordan covers and rounded integer Gram factors
are checked by \path{verify_all.py}.

\subsection{Seven points at the silver ratio}
For four sorted points, if one point beats the sum of the other
three after multiplication by $c$, every larger point does too.
The smallest cannot win when $c<3$.
Thus a four-point subset has a pattern given by the number of
winning targets, from $0$ to $3$.
A pattern on $s$ points assigns such a number to each four-point
subset.

At $c=1+\sqrt2$, joining the realizable patterns level by level gives
$1{,}024$, $3{,}999$ and $1{,}033{,}250$ candidates on five, six and
seven points. Exact Gordan certificates over $\mathbb Q(\sqrt2)$
exclude $972$, $1{,}184$ and $649{,}718$ candidates, leaving
$383{,}532$ seven-point patterns.
The certificate uses the families
$(t,f)=(1,4),(2,4),(3,5),(4,5),(5,6)$, with $59$ blocks in all.
Its matrices have the form
$Q_\sigma=L_\sigma L_\sigma^\top/4^{44}$ for integer $L_\sigma$,
and its exact maximum is the fraction in \eqref{eq:s3-seven}.

An exact verifier checks every retained pattern.
A second verifier, written independently from the definitions and
sharing no code with the first, rebuilds the cover and checks every
exclusion in exact $\mathbb Q(\sqrt2)$ arithmetic.
It evaluates all scores in floating point and the $300$ largest
exactly; the $300$th lies $5\cdot10^{-5}$ below the maximum.
Both routes find the same maximum.
The repository contains the exact verifier and its two inputs in
\path{anc/proofs/s3/silver_ratio_seven_point/}. The script
\path{verify7_k.py} reads the losslessly compressed JSON inputs
directly, rebuilds all three covers, and checks every retained score.
It needs NumPy for batching exact integer arithmetic; no numerical
solver or search checkpoints are required.
The arXiv ancillary package omits the large input
\path{cert7_full_k.json.gz}; reproducing this bound requires obtaining
that file from the repository.

\subsection{The ladder lower bound at the silver ratio}\label{app:s3-lower}
Put $c_*=1+\sqrt2$. One explicit strict ladder, used in the Lean
proof, is
\[
 a_j=\sqrt2-(\sqrt2-1)4^{-j}\quad(0\le j\le N).
\]
Then $a_0=1$, the points increase strictly, and $a_N<\sqrt2$.
The two useful margins are explicit:
\[
 c_*a_{j+1}-(2+a_j)
   =\frac{(\sqrt2-1)(3-\sqrt2)}{4^{j+1}}>0,
 \qquad c_*-(1+a_j)=(\sqrt2-1)4^{-j}>0.
\]
For nonnegative masses $p+q+r=1$ with $p<1$, put
\[
 \nu_N=p\delta_0+q\delta_1+\frac rN\sum_{j=1}^N\delta_{a_j}.
\]
Counting the guaranteed successes as in Section~\ref{sec:three}
gives a base probability at least $F_N(p,q,r)$. This count holds
for every $N\ge1$; it is not inferred from checking finitely many
ladder lengths.

To make the repetition across scales explicit, write $M=a_N$ and
let $d>0$ be the minimum positive value of
$c_*x_4-x_1-x_2-x_3$ on the finite support of $\nu_N$.
Choose
\[
 0<\rho\le\min\left\{\frac{d}{6M},\frac{1}{2c_*M},\frac12\right\}.
\]
Replacing a zero draw by a draw from a $\rho$-scaled copy adds
at most $\rho M$ to each summand. Every old strict success
therefore retains margin at least $d/2$. A positive base summand
is at least $1$, while a lower-scale target satisfies
$c_*\rho M\le1/2$, so it cannot win against such a summand.
Repeating the construction $J$ times gives a finite law with
success probability at least
\[
 F_N(p,q,r)\sum_{k=0}^{J-1}p^{4k}.
\]
Letting $J\to\infty$ proves
$S_3(c_*)\ge F_N(p,q,r)/(1-p^4)$.
The scale ratio may depend on $N$; no common positive ratio for
all ladder lengths is required.

Write $F_\infty=\lim_{N\to\infty}F_N$.
Here is an exact specification of the masses used for $L_3$.
Let $\theta$ be the unique root in
$(0.81558932228969,0.81558932228970)$ of
\[
\begin{aligned}
 P(T)={}&72T^{14}-132T^{13}-526T^{12}+207T^{11}
       +5374T^{10}-8269T^9\\
       &-11547T^8+53262T^7-87309T^6+88506T^5
       -60318T^4\\
       &+27243T^3-7695T^2+1215T-81.
\end{aligned}
\]
Set
\[
 z=\frac{\theta(2\theta^2-6\theta+3)}{3(\theta-1)(3\theta-1)},
 \qquad
 (p,q,r)=\frac{(z,\theta,1-\theta)}{1+z},
 \qquad L_3=\frac{F_\infty(p,q,r)}{1-p^4}.
\]
These masses are positive and sum to one, and
\[
 \frac{F_N(p,q,r)}{1-p^4}
 =L_3-\frac{3q^2r^2}{2N(1-p^4)}\longrightarrow L_3.
\]
Consequently $S_3(c_*)\ge L_3$.
Exact root counting isolates $\theta$; elimination and
irreducibility modulo $37$ show that $L_3$ has degree $14$.
Neither global optimality of the chosen masses within the ladder
family nor sharpness over all laws is needed for this lower bound.

The notes, Python programs, saved rational witnesses, and Lean
proofs are in \path{anc/proofs/s3/silver_ratio_lower/}.
The standard-library script \path{verify_lower.py} enumerates
the endpoint ladders of lengths $1,3,4,5$ in exact
$\mathbb Q(\sqrt2)$ arithmetic and checks the saved rational
witnesses and scale ratios below the endpoint.
The SymPy script \path{verify_field.py} checks the algebraic
identities, root isolation, and degree assertions.
The Lean files \path{LadderGeometry.lean}, \path{Ladder.lean}, and
\path{LowerBound.lean} prove the geometry, count, and limiting
lower bound for all lengths. \path{Amplification.lean} proves
existence of a preserving scale by continuity, using only finite
laws. \path{OptimizedLadder.lean} defines the algebraic masses
and also proves, with separate rational masses, a lower bound
strictly greater than $0.4952063382580596903$.
The final statements in \path{Theorems.lean} concern actual
finite-support probability laws; \path{Audit.lean} lists their
axiom dependencies.

\end{document}